\RequirePackage[leqno]{amsmath}
\documentclass[leqno]{amsart}
\usepackage{color}
\usepackage{graphicx}
\usepackage{subfigure}
\usepackage{amsfonts}
\usepackage{enumitem}
\setlist[itemize]{topsep=0pt,after=\vspace{1.5\baselineskip}}
\usepackage[colorlinks=true]{hyperref}
\usepackage{empheq}
\usepackage[T1]{fontenc}
\usepackage[latin1]{inputenc}

\def\R{\mathbb R} \def\N{\mathbb N} 
\def\Normchi{L^\infty([\gamma,\Gamma])}
\newtheorem{theorem}{Theorem}[section]

\newtheorem{lemma}[theorem]{Lemma}

\newtheorem{remark}{Remark}

\newcommand{\Vast}{\bBigg@{5}}

\title[On the solvability of a general Keller-Segel system] %Use the shortened version of the full title
      {Solvability of a Keller-Segel system with signal-dependent sensitivity and essentially sublinear production}
\author[G. Viglialoro and T.E. Woolley]{}
\keywords{Nonlinear parabolic systems, chemotaxis, global existence, boundedness. \\
\textit{$^\sharp$Corresponding author}: giuseppe.viglialoro@unica.it}
\begin{document}
\maketitle
\maketitle

\centerline{\scshape Giuseppe Viglialoro$^{1,\sharp}$ \and Thomas E. Woolley$^2$}
\medskip
{\footnotesize
 \centerline{$^1$Dipartimento di Matematica e Informatica}
 \centerline{Universit\`{a} di Cagliari}
 \centerline{Viale Merello 92, 09123. Cagliari (Italy)}
 \medskip
 \footnotesize
 \centerline{$^2$Cardiff School of Mathematics}
   \centerline{Cardiff University}
   \centerline{Senghennydd Road,  Cardiff, CF24 4AG (United Kingdom)}
 \medskip
}

\bigskip
\begin{abstract}
In this paper we consider the zero-flux chemotaxis-system 
\begin{equation*}
\begin{cases}
u_{t}=\Delta u-\nabla \cdot (u \chi(v)\nabla v) & \textrm{in}\quad  \Omega\times (0,\infty), \\
0=\Delta v-v+g(u) & \textrm{in}\quad \Omega\times (0,\infty),\\
%\frac{\partial u}{\partial \nu}=\frac{\partial v}{\partial \nu}=0, & x\in \partial \Omega, t>0, \\
%u(x,0)=u_{0}(x) \geq 0 \,\,\textrm{and}\,\,\,v(x,0)=v_{0}(x) \geq 0,& x\in  \Omega,
\end{cases}
\end{equation*}
in a smooth and bounded domain $\Omega$ of $\mathbb{R}^2$. The chemotactic sensitivity $\chi$ is a general nonnegative function from $C^1((0,\infty))$ whilst $g$, the production of the chemical signal $v$, belongs to $C^1([0,\infty))$ and satisfies $\lambda_1\leq g(s)\leq \lambda_2(1+s)^\beta$, for all $s\geq 0$, $0\leq\beta\leq \frac{1}{2}$ and $0<\lambda_1\leq \lambda_2.$ It is established that no chemotactic collapse for the cell distribution $u$ occurs in the sense that any arbitrary nonnegative and sufficiently regular initial data $u(x,0)$ emanates a unique pair of global and uniformly bounded functions  $(u,v)$ which classically solve the corresponding initial-boundary value problem.  Finally, we illustrate the range of dynamics present within the chemotaxis system by means of numerical simulations.
\end{abstract}
% % % % % % % % % % %
\section{Introduction, motivations and main result}
\subsection{Chemotaxis models: general overview and results}\label{Intro}
In this paper we  focus our attention on one of the innumerable variants of the landmark models proposed by Keller and Segel (\cite{K-S-1970} and \cite{Keller-1971-TBC}), which idealize chemotaxis phenomena (largely spread in the physical and biological sciences), describing situations where the motion of a certain individual cells $u=u(x,t)$ is strongly influenced by the presence of a chemical signal $v=v(x,t)$. Precisely a very general mathematical formulation of these models involves two coupled  partial differential equations and reads
\begin{equation}\label{VeryGeneralKeller-Segel}
\left\{ \begin{array}{ll}
u_{t}=  \nabla \cdot (A(u,v)\nabla u-B(u,v)\nabla v)+C(u,v)& \textrm{in}\quad \Omega\times (0,\infty), \\
\tau v_{t}= \Delta v+E(u,v)&\textrm{in}\quad \Omega\times (0,\infty),\\
\frac{\partial u}{\partial \nu}=\frac{\partial v}{\partial \nu}=0 & \textrm{in}\quad \partial \Omega, t>0, \\
u(x,0)=u_{0}(x) \geq 0 \,\,\textrm{and}\,\,\,v(x,0)=v_{0}(x) \geq 0& x\in  \Omega,
\end{array} \right. 
\end{equation}
where $\Omega \subset \R^n$, with $n\geq 1$, is a bounded domain with smooth boundary, $\tau\in\{0,1\}$ and $A, B, C$ and $E$ are sufficiently regular functions of their arguments.  Additionally, $u_0(x)$ and  $v_0(x)$ are the initial cell and chemical distributions, $\partial/\partial \nu$ stands for the outward normal derivative on $\partial \Omega$, so that the zero-flux boundary conditions on both $u$ and $v$ indicate the total insulation of the domain.

Strong numerical methods support real experiments and observations indicating that the aforementioned movement may eventually degenerate to aggregation processes, where an uncontrolled gathering of cells at certain spatial locations is perceived as time evolves: the so called \textit{chemotactic collapse}. In particular, such a coalescence phenomena can be intuitively justified in this way: for $B(u,v)>0$ the first equation in \eqref{VeryGeneralKeller-Segel} indicates that while the cell population naturally diffusing (through the law of $A$) and growing/dying (through the law of $C$), it is additionally driven by the concentration gradient of the chemical signal in the opposite direction of the diffusion (due to the positivity of $B$).

By the mathematical point of view, the chemotactic collapse implies that possibly $u$, in a particular  instant (blow-up time), becomes unbounded in one or more points of its domain. The appearance of this instability may be tied to the size of the initial data $(u_0,v_0)$, the growth rate of the cell distribution $u$ induced by the source term $C(u,v)$, the mutual interplay between the diffusion $A(u,v)$ and the chemotactic sensitivity $B(u,v)$ (also in terms of the space dimension) as well as the production/degradation rate of the chemical signal $v$ given by $E(u,v)$.  

For the linear diffusion case $A(u,v)\equiv 1$, and $B(u,v)=\chi u$, $\chi>0$, in absence of source $C(u,v)$ and with production $E(u,v)=-v+u$, for $\tau=1$ (parabolic-parabolic case), in \cite{OsYagUnidim} it is shown that in one-dimensional domains, all the solutions of  \eqref{VeryGeneralKeller-Segel} are global and uniformly bounded in time, while in the $n$-dimensional context, with $n\geq 2$, unbounded solutions to the same problem have been discovered (see, for instance, \cite{HorstWang} and \cite{WinklerFiniteTeimeBlowUpHigher}). Similarly, for the parabolic-elliptic case ($\tau=0$), in \cite{JaLu} for radial solutions and in \cite{Nagai} for non-radial, the authors prove that for $n=2$ a certain threshold value given by the  product between the chemosensitivity $\chi$ and the initial mass $\int_\Omega u_0$ decides whether the solution can blow up in finite time  or exists for all time $t>0$. 

By virtue of this, let us observe that for positive chemical and cell distributions, the expression $E(u,v)=-v+u$  in \eqref{VeryGeneralKeller-Segel} shows how an increasing of the cells favours a secretion of the signal, which depending on the expression of the chemosensitivity $B(u,v)$ might strongly contrast the smoothing and equilibrating effect of the  diffusion $A(u,v)$. In this sense an abundant literature concerning existence and properties of global, uniformly bounded or  blow-up solutions, is available; for a complete picture, we suggest the largely cited contributions \cite{BellomoEtAl} and \cite{Hillen2009UGP} where, inter alia, reviews of various models about Keller-Segel-type systems are presented and analysed.

Finally, exactly in order to provide a more comprehensive picture of the general situation tied to this balance between the destabilizing and stabilizing effects ($A$ versus $B$), the presence of an absorptive  logistic-type source of the type $C(u,v)\simeq ku-\mu u^\delta$,  for $k \in \R$, $\mu>0$ and $\delta>1$ in \eqref{VeryGeneralKeller-Segel}  may, or may not, have a certain relevance on the dynamics of the system. We mention the papers  \cite{Lankeit,TelloWinkParEl,ViglialoroVeryWeak,ViglialoroBoundnessVeryWeak,ViglialoroWolleyDCDS,WinDespiteLogistic} and \cite{WinklerZAMPLogisticBlowUP} (and references therein), which deal with existence, blow-up and properties of solutions to  \eqref{VeryGeneralKeller-Segel}, for both fully parabolic and parabolic-elliptic versions, and as above for the linear diffusion case $A(u,v)\equiv 1$, and $B(u,v)=\chi u$, $\chi>0$.
\subsection{Some inspiring results related to our research: presentation of the main theorem}
In this investigation the source $C(u,v)$ is taken nil, whereas the sensitivity $B(u,v) = u\chi(v)$ is such that $\chi$ belongs to a general family of positive functions. In particular,
\[B(u,v)=\frac{\chi_0}{v}u,\quad v>0,\quad \textrm{with some}\quad \chi_0>0,\] 
 employed in the so-called Weber-Fechner law expressing the relation between the actual change in the stimulus and the perceived change (\cite{SleemanLevine,OthmerStevens}), is a prototype of such functions, and it presents a singularity at $v=0$ which is the main source of both technical and numerical difficulties. Moreover, and as far as we know this is the novelty of our contribution  (at least in this context), we take the function $E(u,v)$ in such a way that, basically, the corresponding reproduction rate implies a lower increasing of the signal  than that supplied by the linear one $E(u,v)=-v+u$. As a consequence of this, it could be expected that lessening the impact of high values of the cell distribution on the production of the chemical signal may enforce global existence of solutions. To be precise, what could we expect, for instance,  if the chemical were secreted with a sublinear rate depending on the bacterial density? 

Before giving a first answer to this question,  we want to motivate this work starting with an overview of previous achievements regarding some variants of systems like \eqref{VeryGeneralKeller-Segel}, defined in $n$-dimensional domains. In particular, we believe that the coming contributions, all characterised by the presence of singular chemo-sensitivity and expression for $E$ given by $E(u,v)=-v+u$, deserve to be discussed since they put into perspective and also inspire our current investigation: 
\begin{itemize}[leftmargin=*]
\item  for $\tau=1$, $A(u,v)=1$, $B(u,v)=u\chi(v)$, with $\chi(v)=\frac{\chi_0}{v}$, and $C(u,v)=0$, global existence of weak solutions under the assumption $0<\chi_0<\frac{n}{2}$ is proved (\cite{BilerLocalAndGlobal});
\item  for $\tau=0$, $A(u,v)=1$, $B(u,v)=u\chi(v)$, with  $\chi(v)=\frac{\chi_0}{v}$, and $C(u,v)=0$, uniform boundedness and blow-up of radial solutions are positive addressed in \cite{NagaiSenba}; more exactly, solutions are global and remain bounded when either $n\geq 3$ and $0<\chi_0<\frac{n}{n-2}$ or $n=2$ and $\chi_0>0$ is arbitrary, whilst for $n\geq 3$, $0<\chi_0<\frac{n}{n-2}$ and $\int_\Omega u_0|x|$ is sufficiently small, the solution blows up in finite time; 
\item  for $\tau=0$, $A(u,v)=1$, $B(u,v)=u\chi(v)$, for $0<\chi(v)\leq \frac{\chi_0}{v^k}$, with $k\geq 1$ and $v>0$, and $C(u,v)=0$,  global existence and uniform boundedness of classical non-radial solutions are discussed in \cite{FujieWinklerYokotaSignalDependent}, where it is shown that the system possesses a unique global classical solution, uniformly bounded if $0<\chi_0<\frac{2}{n}$ $(k=1)$ and $0<\chi_0<\frac{2}{n}\frac{k^k}{(k-1)^{k-1}}\gamma^{k-1}$ $(k>1)$, $\gamma>0$ being a constant depending on the data.
\item   for $\tau=0$, $A(u,v)=1$, $B(u,v)=u\chi(v)$, with  $\chi(v)=\frac{\chi}{v}$, and $C(u,v)=ru-\mu u^2$, $r\in\R$, $\chi,\mu>0$, it is proved that in two-dimensional domains the logistic kinetics ensures global existence of classical solutions even for arbitrary large $\chi$ and any $\mu>0$. Additionally, it is shown that if $r$ is larger with respect to some expression of $\chi$ such solutions are also bounded (\cite{FujieWinklerTomomiParaEllip}); 
\item   for $\tau=1$, $A(u,v)=1$, $B(u,v)=u\chi(v)$, with $\chi(v)=\frac{\chi}{v}$, and $C(u,v)=0$, in \cite{LankeitANewApproach} uniform boundedness of global classical solutions is shown in the two-dimensional setting and for $\chi\in(0,\chi_0)$, for some $\chi_0>1$. 
\item  for $\tau=0$, $A(u,v)=(u+1)^{m-1}$, $B(u,v)=u(u+1)^{\alpha-1}\chi(v)$, with $0<\chi(v)\leq \frac{\chi_0}{v^k}$, $k\geq 1$ and $v>0$, and $C(u,v)=0$, for $m,\alpha \in \R$ such that $\alpha \leq \max\{m,\frac{m+1}{2}\}$, global existence and uniform boundedness of classical solutions are proved, provided some smallness assumptions on $\chi_0$ are satisfied (\cite{ViglialoroPreprintNonLinear}).
\end{itemize}  
In accordance with these premises, we intend to enhance the knowledge of the mathematical analysis of general chemotaxis-systems, by studying this problem
\begin{equation}\label{problem}
\begin{cases}
u_{ t}=\Delta u-\nabla \cdot (u\chi(v)\nabla v) & \textrm{in}\quad \Omega\times (0,\infty), \\
0=\Delta v-v+g(u) & \textrm{in}\quad \Omega\times (0,\infty),\\
\frac{\partial u}{\partial \nu}=\frac{\partial v}{\partial \nu}=0 &\textrm{on}\quad \partial \Omega\times (0,\infty), \\
u(x,0)=u_{0}(x) \geq 0& x\in  \Omega,
\end{cases}
\end{equation}
where $\Omega$ is a smooth and bounded domain of $\mathbb{R}^2$ and $u_0(x)=u(x,0)$ is a nonnegative function such that $0\not \equiv u_0\in C^0(\bar{\Omega})$. 
% % % % % % % % % % % %
%\subsection{Presentation of the main result and plan of the paper}\label{IntroductionSection} 
Additionally, $0<\chi\in C^1((0,\infty))$ and $g\in C^1([0,\infty))$ satisfies this \textit{essentially sublinear growth} (see Remark \ref{RemarkEssentiallySublinear} below):
\begin{equation}\label{AssumptionsOng}
 0<\lambda_1 \leq g(s)\leq \lambda_2(1+s)^\beta \quad \textrm{for} \quad s\geq 0, \quad 0\leq\beta\leq \frac{1}{2}  \quad \textrm{and}\quad 0<\lambda_1\leq \lambda_2.
\end{equation}
We prove the existence and uniqueness of a global and bounded classical solution to problem \eqref{problem}, and precisely we show that  the high sublinear action induced on $v$ by $g$ exerts a certain smoothing effect on $u$ and it is sufficient to prevent  $\delta$-singularities formation for the cell distribution $u$, even for widely large initial cells' density $u_0(x)$ or strong sensitivity effects.

This conclusion is mathematically formulated in our main theorem:
\begin{theorem}\label{MainTheorem}
%{\textcolor{blue}{
Let $\Omega$ be a smooth and bounded domain of $\mathbb{R}^2$, $0<\chi\in C^1((0,\infty))$ and $g\in C^1([0,\infty))$ a function satisfying  \eqref{AssumptionsOng}. Then for any nonnegative initial data $0\not \equiv u_0\in C^0(\bar{\Omega})$, problem \eqref{problem} admits a unique global classical solution $(u, v)$. Moreover, both $u$ and $v$ are bounded in $\Omega \times (0,\infty).$
\end{theorem}
\begin{remark}\label{RemarkEssentiallySublinear}
For the chemotaxis model \eqref{problem} with linear production, i.e. $g(u)=u$, it is seen from the conservation of the total cell mass (see Lemma \ref{LocalExistenceLemma}), that is $\int_\Omega u=\int_\Omega u_0$, that also the same property for the chemical $v$ holds; indeed, by integrating over $\Omega$ the second equation we have $\lVert v\Vert_{L^1(\Omega)}=\int_\Omega g(u)=\int_\Omega u=\int_\Omega u_0>0$ throughout the time. Conversely,  for chemotaxis models with sublinear production, $g(u)=u^\vartheta$ with $0<\vartheta<1$,  this is no longer true. In fact, Hölder's inequality implies $\lVert v\Vert_{L^1(\Omega)}=\int_\Omega g(u)=\int_\Omega u^\vartheta\leq (\int_\Omega u_0)^\vartheta |\Omega|^{1-\vartheta}$, which does not exclude the possibility of vanishing for $\int_\Omega v$ at some time. Since in order to save the chemosensitivity $\chi=\chi(v)$ from  singularities  we have to avoid this last scenario (see Lemma \ref{LemmaLowerBoundv}),  we assume   $g(s)\geq \lambda_1>0$, for $s>0$. Subsequently, for the production source $g$ restricted to grow no faster than $u^\beta$ at infinity (essentially sublinear growth), i.e. $g(s) \leq \lambda_2(1+ s)^\beta$ for all $s>0$ and some $\lambda_2>0$,  the assumption $\lambda_2\geq \lambda_1>0$ makes consistent the lower and upper bounds in \eqref{AssumptionsOng}.
\end{remark}
The rest of the paper is structured as follows. First, in $\S$\ref{PreliminariesSection}, we collect some necessary and preparatory material, then, in $\S$\ref{ExistenceSolutionRegularizingSection}, we prove the local existence and uniqueness of a classical solution to \eqref{problem} and some of its properties. Successively, in $\S$\ref{FromLocalToGlobaSection}, we establish how to ensure globability and boundedness of local solutions using their $L^p$-boundedness. Such a bound is derived in $\S$\ref{EstimatesAndProofSection}, which represents the main part of this report and that concludes with the proof of Theorem \ref{MainTheorem}. Finally, the theoretical results presented here are investigated numerically in $\S$\ref{SimuationsSection}, where simulations are used to detect critical exponents for $\beta$ which
delineate regions where different asymptotic behaviours of solutions to the same system \eqref{problem} may manifest. 
% % % % % % % % % % % % % % % % % % % % % % % % % % % % % % % % % % % % % % % % % % % % % % % %
\section{Preliminaries and auxiliary tools}\label{PreliminariesSection}
The coming results are supportive in the proof of the main theorem of this paper. To be precise, we mainly summarize and derive some general functional inequalities, also tied to elliptic regularity theory.

Let us first recall a special case of the well-known Gagliardo-Nirenberg inequality which will be used through the paper to prove the main theorem. 
\begin{lemma} (Gagliardo-Nirenberg inequality)\label{InequalityG-NLemma}
Let $\Omega$ be a smooth and bounded domain of $\R^2$. Then there is a constant $C_{GN}>0$ such that the following inequality holds: 
With $\mathfrak{q},\mathfrak{s}\in[1,2]$, $\mathfrak{p}\in [2,4]$ and $\theta=1-\frac{\mathfrak{q}}{\mathfrak{p}}\in[0,1)$,  
\begin{equation}\label{InequalityTipoG-N} 
\| f \|_{L^{\mathfrak{p}}(\Omega)} \leq C_{GN} ( \| \nabla f \|_{L^{2}(\Omega)}^{\theta} \| f \|_{L^{\mathfrak{q}}(\Omega)}^{1 - \theta}+  \| f \|_{L^{\mathfrak{s}}(\Omega)})
\end{equation}
is satisfied for all $f\in L^\mathfrak{q}(\Omega)$ with $\nabla f\in L^2(\Omega)$,  
\begin{proof}
See \cite[p. 126]{NirenbergOnEllipticPDE}.
\end{proof}
\end{lemma}
The following lemma is fundamental in our computations and its validity is restricted to two-dimensional settings, which are those where our main problem is studied. Its proof is, fundamentally, a reformulation of \cite[Lemmas 4.3. and 4.4.]{FujieWinklerTomomiParaEllip} which we adapted to our presentation in order to make the present article more self-contained.
\begin{lemma}\label{LemmaTwoDimensional}
Let $\Omega$ be a smooth and bounded domain of $\R^2$.  Then there exists a positive constant $\hat{C}$ such that for all $p\in(1,2)$ and $f\in C^2(\bar{\Omega})$, with $\frac{\partial f}{\partial \nu}=0$ on $\partial \Omega$, holds that
\begin{equation*} %\label{EstimateNablaf2pplus2} 
\lVert \nabla f \rVert_{L^{2p+2}(\Omega)}\leq \hat{C}\lVert -\Delta f + f\rVert^{\frac{1}{2}}_{L^{p+1}(\Omega)}\lVert f \rVert^{\frac{1}{2}}_{L^2(\Omega)}.
\end{equation*}
\begin{proof}
Given $\psi\in C^2(\bar{\Omega})$, we apply \eqref{InequalityTipoG-N} with $\mathfrak{p}=4$ and $\mathfrak{q}=2$ to obtain
\begin{equation}\label{InequalityG-N_TwoDimension} 
\| \psi  \|_{L^{4}(\Omega)}^4 \leq c_1  \| \nabla \psi  \|_{L^{2}(\Omega)}^{2} \| \psi \|_{L^{2}(\Omega)}^{2}+  c_1\| \psi \|_{L^{2}(\Omega)}^4,
\end{equation}
where $c_1=(2C_{GN})^4$, having also made use of  
\begin{equation}\label{AlgebraicInequality2toalpha}  
(a+b)^\alpha\leq 2^\alpha(a^\alpha+b^\alpha)\quad \textrm{for any} \quad a,b\geq 0, \alpha>0.
\end{equation}
In addition, for any $\varphi \in C^1(\bar{\Omega})$ we have that  for all $p\in (1,2)$ the function $\psi=|\varphi|^\frac{p+1}{2}$ belongs to $W^{1,2}(\Omega)$ and is such that $|\nabla \psi|=\frac{p+1}{2}|\varphi|^\frac{p-1}{2}|\nabla \varphi|;$ consequently, inequality \eqref{InequalityG-N_TwoDimension}  explicitly reads
\begin{equation}\label{InequalityFirstForTwoDimensional}
\int_\Omega |\varphi|^{2(p+1)}\leq  c_1\frac{(p+1)^2}{4}\bigg(\int_\Omega |\varphi|^{p-1}|\nabla \varphi|^2\bigg)\bigg(\int_\Omega |\varphi|^{p+1}\bigg)+c_1\bigg(\int_\Omega |\varphi|^{p+1}\bigg)^2.
\end{equation}
Now, the Hölder inequality enables us to get
\begin{equation*} %\label{InequalityFirstHalfForTwoDimensional}
\int_\Omega |\varphi|^{p-1}|\nabla \varphi|^2 \leq \bigg(\int_\Omega |\nabla \varphi|^{p+1}\bigg)^\frac{2}{p+1}\bigg(\int_\Omega |\varphi|^{p+1}\bigg)^\frac{p-1}{p+1},
\end{equation*}
and
\begin{equation*} %\label{InequalitySecondForTwoDimensional}
\bigg(\int_\Omega |\varphi|^{p+1}\bigg)^\frac{2p}{p+1} =\bigg(\int_\Omega |\varphi|^\frac{p^2-1}{p}|\varphi|^\frac{p+1}{p}\bigg)^\frac{2p}{p+1}\leq \bigg(\int_\Omega | \varphi|^{2(p+1)}\bigg)^\frac{p-1}{p+1}\int_\Omega |\varphi|^{2}.
\end{equation*}
Further, by inserting these last two relations into \eqref{InequalityFirstForTwoDimensional}, a proper decomposition and the inequality $(p+1)^2>4$, valid for all $1<p<2$, infer
\begin{equation*} %\label{InequalityFirstForTwoDimensionalBis}
\begin{split}
\int_\Omega |\varphi|^{2(p+1)}&\leq  c_1\frac{(p+1)^2}{4}\bigg[\bigg(\int_\Omega |\nabla \varphi|^{p+1}\bigg)^\frac{2}{p+1}+\bigg(\int_\Omega |\varphi|^{p+1}\bigg)^\frac{2}{p+1}\bigg]  \\ &
\quad \times \bigg(\int_\Omega |\varphi|^{p+1}\bigg)^\frac{2p}{p+1},
\end{split}
\end{equation*}
so that algebraic manipulations yield 
\begin{equation}\label{InequalityThirdForTwoDimensional}
\begin{split}
\lVert \varphi \Vert_{L^{2(p+1)}(\Omega)}^4&=\bigg(\int_\Omega |\varphi|^{2(p+1)}\bigg)^\frac{2}{p+1} =\bigg(\int_\Omega |\varphi|^{2(p+1)}\bigg)^{1-\frac{p-1}{p+1}} \\ &\leq  c_1\frac{(p+1)^2}{4}\bigg[\bigg(\int_\Omega |\varphi|^{p+1}\bigg)^\frac{2}{p+1}
 +\bigg(\int_\Omega |\nabla \varphi|^{p+1}\bigg)^\frac{2}{p+1}\bigg]\int_\Omega |\varphi|^{2}\\ &
\leq c_2\bigg[\bigg(\int_\Omega |\varphi|^{p+1}+\int_\Omega |\nabla \varphi|^{p+1}\bigg)\bigg]^\frac{2}{p+1}\int_\Omega |\varphi|^{2}\\ &
=c_2\lVert \varphi \Vert_{W^{1,p+1}(\Omega)}^2\lVert \varphi \Vert_{L^{2}(\Omega)}^{2},
\end{split}
\end{equation}
where we used, in the last step, $a^\gamma+b^\gamma\leq 2^{1-\gamma}(a+b)^\gamma$, for all $a,b\geq 0$ and $0<\gamma\leq 1$  and set $c_2= c_1\frac{(p+1)^2}{4}2^\frac{p-1}{p+1}$.

On the other hand, let $C>0$ the constant from \cite[(4.6) of Lemmas 4.3.]{FujieWinklerTomomiParaEllip}; for $f$ as in our hypothesis and $p\in(1,2)$ we can rely on such inequality and obtain
\begin{equation}\label{PrimeraDesigualdadInNorm}
\lVert f \rVert_{W^{2,p+1}(\Omega)}\leq C \lVert -\Delta f + f \rVert_{L^{p+1}(\Omega)}.
\end{equation}
In addition, the same properties of $f$ also allow us to set $\varphi=|\nabla f|$ in \eqref{InequalityThirdForTwoDimensional} so to have for $c_3=c_2^\frac{1}{4}$
\begin{equation}\label{SegundaDeisgualdadInNorm} 
\lVert \nabla f \rVert_{L^{2(p+1)}(\Omega)}\leq c_3 \lVert f \rVert_{W^{2,p+1}(\Omega)}^\frac{1}{2}\lVert \nabla f \rVert_{L^{2}(\Omega)}^\frac{1}{2},
\end{equation}
where we considered \eqref{AlgebraicInequality2toalpha} and the fact that 
\[
\begin{split}
\lVert \nabla f\Vert_{W^{1,p+1}(\Omega)}^{p+1}&=
\int_\Omega |\nabla f|^{p+1}+\int_\Omega|\nabla  |\nabla f||^{p+1}\\ & \leq \int_\Omega |f|+\int_\Omega |\nabla f|^{p+1}+\int_\Omega|\nabla  |\nabla f||^{p+1}\leq \lVert  f\Vert_{W^{2,p+1}(\Omega)}^{p+1}.
%uu
\end{split}
\]
We conclude the proof by combining \eqref{PrimeraDesigualdadInNorm} and \eqref{SegundaDeisgualdadInNorm}, being $\hat{C}=c_3\sqrt{C}.$
\end{proof}
\end{lemma}
Conforming to the comments in Remark \ref{RemarkEssentiallySublinear}, in the next result we will establish a lower bound for the second component of solutions to the parabolic-elliptic Keller-Segel system \eqref{problem}. In particular we derive a quantitative estimate on positivity of solutions to the Neumann problem for the Helmholtz equation with nonnegative inhomogeneity having given norm in $L^1(\Omega)$.
\begin{lemma}\label{LemmaLowerBoundv}
For any $n\in \N$, let $\Omega$ be a smooth and bounded domain of $\mathbb{R}^n$ and $0 \not\equiv w\in C^0(\bar{\Omega})$ a nonnegative function.  If $f\in C ^2(\bar{\Omega})$ is a solution of 
\[
\begin{cases}
-\Delta f + f=w & \textrm{in}\quad  \Omega,\\
\frac{\partial f}{\partial \nu}=0& \textrm{on}\quad \partial \Omega,
\end{cases}
\]
then there exists some positive constant $\eta$ such that 
\[f\geq \eta \int_\Omega w>0\quad \textrm{in}\quad \Omega.\]
\begin{proof}
The proof is a consequence of the the positivity of the Green function to the Helmholtz equation.
\end{proof}
\end{lemma}
We will also make use of the following elementary proof:
\begin{lemma}\label{BoundsInequality}
Let $y$ be a positive real number verifying $ y \leq k(y^l+1)$ for some $k>0$ and $0<l<1$. Then $y\leq \max\{1,(2k)^\frac{1}{1-l}\}.$
\begin{proof}
 Since for $y\leq 1$ there is nothing left to show, we suppose $y\geq 1$. Then $y^l \geq 1$ so that $y \leq k(y^l+1) \leq 2k y^l$ and hence $y\leq  (2k)^\frac{1}{1-l}$.
\end{proof}
\end{lemma}
\section{Existence of local-in-time solutions and main properties}\label{ExistenceSolutionRegularizingSection}
We open this section with a lemma concerning the local-in-time existence of classical solutions $(u,v)$ to system \eqref{problem}. The proof is developed by adapting well-established methods involving an appropriate fixed point framework and standard parabolic and elliptic regularity results. Through the same lemma we also are able to achieve an important conservation of mass property and, relying on Lemma \ref{LemmaLowerBoundv}, a crucial \textit{uniform-in-time} estimate for the component $v$ of the solution, which ensures the uniform bound of \textit{any} signal-dependent sensitivity $\chi=\chi(v)$ taken from $C^1((0,\infty)).$
\begin{lemma}\label{LocalExistenceLemma}  
For any $n\in \N$, let $\Omega$ be a smooth and bounded domain of $\mathbb{R}^n$, $0<\chi\in C^1((0,\infty))$ and and $g\in C^1([0,\infty))$ a function satisfying  \eqref{AssumptionsOng}. Then for any nonnegative initial data $0\not \equiv u_0\in C^0(\bar{\Omega})$, problem \eqref{problem} admits a unique local-in-time classical solution 
\begin{equation*}
(u,v)\in (C^0([0,T_{max});C^0(\Omega))\cap C^{2,1}(\bar{\Omega}\times (0,T_{max})))^2,
\end{equation*}
where $T_{max}\in(0,\infty]$, denoting the maximal existence time, is such that if $T_{max}<\infty$ necessarily
\begin{equation}\label{extensibility_criterion_Eq} 
\limsup_{t\nearrow T_{max}}\lVert u (\cdot,t)\rVert_{L^\infty(\Omega)}=\infty.
\end{equation}
Moreover, for some $\gamma,\Gamma>0$ we have for all $(x,t)$ in $\Omega \times (0,T_{max})$
\begin{equation}\label{MaximumPricnipleRelation} 
u\geq 0, \quad \gamma \leq v \leq \Gamma  \quad \textrm{and}\quad  \lVert \chi(v) \rVert_{L^\infty(\Omega \times (0,T_{max}))}\leq \lVert \chi \rVert_{L^\infty([\gamma,\Gamma])},
\end{equation}
and also
\begin{equation}\label{Bound_of_u} 
\int_\Omega u (\cdot ,t)  = m=\int_\Omega u_0>0\quad \textrm{for all}\quad  t\in (0,T_{max}).
\end{equation}
\begin{proof}
\textit{Existence.} For any $T\in (0,1)$, $0\not \equiv u_0\in C^0(\bar{\Omega})$ nonnegative,  and $R:=\lVert u_0\rVert_{L^\infty(\Omega)}+1$, let us consider the Banach space $X:=C^0(\bar{\Omega} \times [0,T])$ and its closed subset
\[S:=\Bigg\{0\leq u \in X\; \bigg \lvert \lVert u (\cdot, t)\rVert_{L^\infty(\Omega)}\leq R \quad \textrm{for all}\quad  t\in[0,T]\\
\Bigg\}.
\]  
For $\hat{u} \in S$, let $v$ be the solution of 
\begin{equation}\label{LocalExistenceProblemforv}
\begin{cases}
-\Delta v+v=g(\hat{u}) &  \textrm{in}\quad \Omega \times (0,T),\\
\frac{\partial v}{\partial \nu}=0 & \textrm{on}\quad \partial \Omega \times (0,T), 
\end{cases}
\end{equation}
and, in turn, let $u$ be the solution of 
\begin{equation}\label{LocalExistenceProblemforu}
\begin{cases}
u_{t}-\Delta u=\nabla \cdot (u \chi(v)\nabla v) &\textrm{in}\quad \Omega \times (0,T), \\
\frac{\partial u}{\partial \nu}=0 &\textrm{in}\quad \partial \Omega \times (0,T), \\
u(x,0)=u_{0}(x) \geq 0 & x\in \Omega.
\end{cases}
\end{equation}
In agreement with these statements, we shall show that for appropriate small $T$,  $\Phi: S\rightarrow S$ defined by $\Phi(\hat{u})=u$ is a compact map such that $\Phi(S)\subset S$; subsequently, due to the convexity of $S$, the Schauder fixed point theorem ensures the existence of  $u\in S$ such that $\Phi(u)=u$.

First, we observe that for a certain fixed $\hat{u} \in S$ well known elliptic regularity results, in conjunction with Morrey's theorem (\cite{BrezisBook}), infer a unique solution $v(\cdot,t)$ to problem \eqref{LocalExistenceProblemforv} in the space $C^{1,\delta}(\Omega)$, for all $\delta \in (0,1)$; this, in particular, implies that  $\nabla v \in L^\infty(\Omega)$ for all $t\in(0,T)$. Again continuing on the property of the solution $v$, because $\hat{u}\in S$, $\hat{u}$ is nonnegative so that from \eqref{AssumptionsOng} we have that $g(\hat{u})$ is well defined and moreover $g(\hat{u})\geq \lambda_1>0$. In this way, an application of Lemma \ref{LemmaLowerBoundv} to  problem \eqref{LocalExistenceProblemforv}, together with the definition of $S$, leads to
\begin{equation}\label{LowerBoundvComponent}
v(x, t) \geq\eta \cdot \int_\Omega g(\hat{u})=\gamma:=\eta \lambda_1|\Omega|>0\quad (x,t)\in\Omega \times (0,T).
\end{equation}
Moreover, besides $v(x,t)\geq \gamma$, the elliptic maximum principle, and again \eqref{AssumptionsOng} and the definition of $S$, provide
\begin{equation}\label{UpperBoundvComponent}
v(x,t)\leq \sup_{\Omega} g(\hat{u})\leq \Gamma :=\lambda_2(1+ R)^\beta \quad (x,t)\in\Omega \times (0,T),
\end{equation}
so that  on account of $\chi\in C^1((0,\infty))$ and $\gamma\leq s:=v(x,t) \leq \Gamma $, with $(x,t)\in\Omega\times (0,T)$, $\chi(s)$ is also from $L^\infty([\gamma,\Gamma])$ and in the specific we have
\begin{equation}\label{BoundednessOfChi}
 \lVert \chi(v) \rVert_{L^\infty(\Omega \times (0,T))}\leq \lVert \chi \rVert_{L^\infty([\gamma,\Gamma])}.
\end{equation}
Subsequently, for some positive constant $c$ (which until the end of this proof might change line by line), using $\hat{u}\in S$ and the gained bounds for $\nabla v$ and $\chi(v)$, \cite[Theorem V 1.1.]{LSUBookInequality} applied to problem \eqref{LocalExistenceProblemforu} implies that $u\in C^{\delta_1,\frac{\delta_1}{2}}(\Omega \times (0,T))$, for some $\delta_1 \in (0,1)$. Hence, 
\[|u(x_1, t_1)-u(x_2,t_2)|\leq c (|x_1-x_2|^{\delta_1}+|t_1-t_2|^\frac{\delta_1}{2})\quad \textrm{for all}\quad x_1,x_2\in \Omega, \, t_1,t_2\in (0,T),\]
that is
\[u(\cdot, t)\leq u_0(\cdot)+ct^\frac{\delta_1}{2}\quad \textrm{for all}\quad \, t\in (0,T).\]
Thereafter
\begin{equation*}
\max_{t\in[0,T]} \lVert u(\cdot,t)\rVert_{L^\infty(\Omega)}\leq \lVert u_0\rVert_{L^\infty(\Omega)}+cT^\frac{\delta_1}{2},
\end{equation*}
and subsequently for $T< c^\frac{-\delta_1}{2}$ we also deduce that 
\begin{equation*}
 \lVert u(\cdot,t)\rVert_{L^\infty(\Omega)}\leq \lVert u_0\rVert_{L^\infty(\Omega)}+1=R \quad \textrm{for all}\quad \, t\in (0,T).
\end{equation*}
Additionally, $\underline{u}\equiv 0$ is a subsolution of the first equation in \eqref{LocalExistenceProblemforu} so that the parabolic comparison principle warrants the nonnegativity of $u$; hence $\Phi $ maps $S$ into itself, compactly since $ C^{\delta_1,\frac{\delta_1}{2}}(\Omega \times (0,T)) \hookrightarrow X$. Let $u$ be a fixed point of $\Phi$; by employing the elliptic and parabolic regularity theory to problems  \eqref{LocalExistenceProblemforv}  and \eqref{LocalExistenceProblemforu} (explicitly \cite[Theorem 9.33]{BrezisBook} and \cite[Theorem V 6.1.]{LSUBookInequality}, respectively), we have $v\in C^{2+\delta_1,\frac{\delta_1}{2}}(\bar{\Omega},[\tau,T])$ and hence $u\in C^{2+\delta_1,1+\frac{\delta_1}{2}}(\bar{\Omega} \times [\tau,T])$, for any $\tau \in (0,T)$. Moreover, by standard bootstrap arguments the solution may be prolonged in the interval  $[0, T_{max})$, with $T_{max}\leq \infty$, $T_{max}$ being finite if and only if \eqref{extensibility_criterion_Eq} holds. In this way, the gained nonnegativity of $u(\cdot,t)$ in $[0,T]$, lower and upper estimates for $v(\cdot,t)$ in $[0,T]$ (see  \eqref{LowerBoundvComponent} and \eqref{UpperBoundvComponent}) and bound for $\chi(v(\cdot,t))$ in $[0,T]$ (see \eqref{BoundednessOfChi}) remain preserved up to $T_{max}$, exactly as claimed in \eqref{MaximumPricnipleRelation}. Finally, an integration over $\Omega$ in the first equation of \eqref{LocalExistenceProblemforu} and the no-flux boundary conditions on both $u$ and $v$ give $\int_\Omega u(\cdot,t)=\int_ \Omega u_0$ for all $t\in (0,T_{max})$, and \eqref{Bound_of_u} is also justified.

\textit{Uniqueness.}  By absurdity let  $(u_1,v_1)$ and $(u_2,v_2)$ be two nonnegative different classical solutions of \eqref{problem} in $\Omega \times (0,T_{max})$ with the same initial data $u_1(\cdot,0)=u_2(\cdot,0)$. In such circumstances, using the equation for $u$ in \eqref{problem}, we have
\begin{equation}\label{FirstInequalityForUniqueness}
\begin{split}
\frac{1}{2}\frac{d}{d t}\int_\Omega (u_1-u_2)^2&+\int_\Omega |\nabla (u_1-u_2)|^2\\ &
=\int_\Omega (u_1\chi(v_1)\nabla v_1-u_2\chi(v_2)\nabla v_2)\cdot \nabla (u_1-u_2)\\ &
\leq \frac{1}{2}\int_\Omega |(u_1\chi(v_1)\nabla v_1-u_2\chi(v_2)\nabla v_2)|^2\\ &
\quad +\frac{1}{2}\int_\Omega |\nabla (u_1-u_2)|^2 \quad \textrm{ for all}\quad t\in(0,T_{max}), 
\end{split}
\end{equation}
and of course 
\begin{equation}\label{InitialDataNilForF}
\int_\Omega (u_1(\cdot,0)-u_2(\cdot,0))^2=0. 
\end{equation}
Now, for all $t\in (0,T_0)$ with $T_0<T_{max}$ we set
\begin{align*}
&s_1=s_1(T_0):=\min\{\lVert u_1 \rVert_{L^\infty(\Omega \times (0,T_0))},\lVert u_2 \rVert_{L^\infty(\Omega \times (0,T_0))}\}, \\ 
&s_2=s_2(T_0):=\max\{\lVert u_1 \rVert_{L^\infty(\Omega \times (0,T_0))},\lVert u_2 \rVert_{L^\infty(\Omega \times (0,T_0))}\},
\end{align*}
and the Mean Value Theorem applied to the function $s\mapsto g(s)$ in the interval $[s_1,s_2]$ infers $g(s_1)-g(s_2)=g'(\bar{s})(s_1-s_2)$ for some $\bar{s} \in (s_1,s_2)$. 
In light of this, through Young's inequality, the second equation of \eqref{problem} provides some positive $C_1$ depending on $T_0$ such that on $(0,T_0)$
\begin{equation}\label{FirstStepToEnsureUniqueness}
\begin{split}
\int_\Omega \lvert \nabla (v_1-v_2)\rvert^2&=-\int_\Omega (v_1-v_2)^2+\int_\Omega (g(u_1)-g(u_2))(v_1-v_2) \\ &
=-\int_\Omega (v_1-v_2)^2+C_1\int_\Omega (u_1-u_2)(v_1-v_2) \\&
\leq \frac{C_1^2}{2}\int_\Omega (u_1-u_2)^2-\frac{1}{2}\int_\Omega (v_1-v_2)^2. %\leq \frac{C_1^2}{2}\int_\Omega (u_1-u_2)^2.
\end{split}
\end{equation}
Additionally, the same elliptic and parabolic regularity results previously used allow us to find some $C_2=C_2(T_0)>0$ such that $|\nabla v_1|\leq C_2$, $|\nabla v_2|\leq C_2$ and $u_2\leq C_2$ on $\Omega \times (0,T_0)$; subsequently, by virtue of the boundedness of $\chi$ in \eqref{MaximumPricnipleRelation} and Hölder's inequality,  some manipulations lead to
\begin{equation*}
\begin{split}
\frac{1}{2}  \int_\Omega & \lvert u_1\chi(v_1)\nabla v_1-u_2\chi(v_2)\nabla v_2 \rvert^2\leq\frac{3}{2} \int_\Omega \lvert (u_1-u_2)^2\chi(v_1)^2|\nabla v_1 |^2\\ &
\quad +\frac{3}{2} \int_\Omega u_2^2\chi(v_2)^2\lvert \nabla (v_1-v_2) \rvert^2 +
\frac{3}{2} \int_\Omega u_2^2(\chi(v_1)-\chi(v_2))^2(v_1-v_2)^2\lvert \nabla (v_2) \rvert^2\\ &
%\leq 2 \lVert \chi \rVert_{L^\infty([\gamma,R])}^2\int_\Omega \lvert u_1-u_2\lvert^2 \rvert\nabla v_1 \rvert^2\int_\Omega \lvert \nabla (v_1-v_2) \rvert^2\\ &
%\quad + 2\lVert \chi' \rVert_{L^\infty([\gamma,R])}^2\int_\Omega \lvert u_2\nabla (v_1-v_2) \rvert^2\int_\Omega \lvert \nabla (v_1-v_2) \rvert^2\\ &
\leq \frac{3}{2}  \lVert \chi \rVert_{L^\infty([\gamma,R])}^2C_2^2 \int_\Omega (u_1-u_2)^2+\frac{3}{2} \lVert \chi \rVert_{L^\infty([\gamma,R])}^2C_2^2\int_\Omega \lvert \nabla (v_1-v_2)\rvert^2\\ & 
\quad +\frac{3}{2}  C_2^4\lVert \chi' \rVert_{L^\infty([\gamma,R])}^2\int_\Omega  (v_1-v_2)^2\quad \textrm{for all}\quad t\in (0,T_0),
\end{split}
\end{equation*} 
where we also applied the Mean Value theorem to the function $s\mapsto\chi(s)$, with $s\in[s_1,s_2]$, and used the  inequality  $(a+b+c)^2\leq  3 (a^2+b^2+c^2)$, valid for all $a,b,c \in \R.$
Hence, we can write for all $t\in(0,T_{max})$ 
\begin{equation}\label{UniquenessFromChi}
\begin{split}
\frac{1}{2}\int_\Omega &|(u_1\chi(v_1)\nabla v_1-u_2\chi(v_2)\nabla v_2)|^2\leq \\&
\quad C_3 \bigg(\frac{1}{2}\int_\Omega (u_1-u_2)^2+\int_\Omega \lvert \nabla (v_1-v_2)\rvert^2 +\frac{1}{2}\int_\Omega  (v_1-v_2)^2\bigg),
\end{split}
\end{equation}
where $C_3=C_3(T_0)= 3C_2^2\max\{\lVert \chi \rVert_{L^\infty([\gamma,R])}^2,\lVert \chi' \rVert_{L^\infty([\gamma,R])}^2C_2^2\}$.
%\begin{equation}\label{UniquenessFromChi}  
%\begin{split}
%I_1 &\leq\frac{ \lVert \chi \rVert_{L^\infty([\gamma,R])}C_2^2}{2}\int_\Omega (u_1-u_2)^2\\ &
%\quad +\left(\frac{ \lVert \chi \rVert_{L^\infty([\gamma,R])}C_2^2}{2}+C_2 \lVert \chi \rVert_{L^\infty([\gamma,R])}\right)\int_\Omega  \lvert \nabla (v_1-v_2)\rvert^2.
%\end{split}
%\end{equation}
%Additionally, application of the Mean Value Theorem, the Young inequality and the boundedness of $u_1$ and $u_2$ in $\Omega \times (0,T_0)$ provide some $C_3=C_3(T_0)>0$
%\begin{equation}\label{UniquenessFromLogistic} 
% \int_\Omega (f(u_1)-f(u_2)) (v_1-v_2)\leq \frac{C_3\rho_2}{2} \int_\Omega (u_1-u_2)^2+\frac{C_3}{2\rho_2}\int_\Omega (v_1-v_2)^2,
%\end{equation}
%with some $\rho_2>0$. 

Finally, coming back to \eqref{FirstInequalityForUniqueness}, and plugging in it \eqref{UniquenessFromChi}  and \eqref{FirstStepToEnsureUniqueness} we arrive at
%\begin{equation*}
%C_3=C_3(T_0)=\max\bigg\{C_1,2 \lVert \chi \rVert_{L^\infty([\gamma,R])} C_2\bigg(1+\frac{C_2^4}{2C_1}\lVert \chi \rVert_{L^\infty([\gamma,R])}\bigg)\bigg\}
%\end{equation*}
this initial problem 
\begin{equation}\label{F_eqn}
\frac{d}{d t} \mathcal{F}\leq (C_3+C_3C_1^2) \mathcal{F} \quad t\in (0,T_0), \quad \mathcal{F}(0)=0,
\end{equation}
where $\mathcal{F}(t):=\int_\Omega (u_1-u_2)^2$ and where \eqref{InitialDataNilForF} has been also considered. Since \eqref{F_eqn} admits the unique solution $ \mathcal{F}\equiv 0$ on $(0,T_0)$, due to the arbitrary of $T_0$, we attain $u_1=u_2$ on $(0,T_{max})$ and hence, by using again \eqref{FirstStepToEnsureUniqueness}, also $v_1=v_2$  on $(0,T_{max})$.
\end{proof}
\end{lemma}
\section{From local to global-in-time and bounded solutions}\label{FromLocalToGlobaSection}
The forthcoming important result shows how to achieve  uniform-in-time boundedness of solutions from their $L^p$-boundedness, for some suitable $p>1$. As a consequence, in order to establish our main theorem, it will be successively sufficient the derivation of such a $L^p$ bound. 
\begin{lemma}\label{FromLocalToGLobalBoundedLemma}
Under the assumptions of Lemma \ref{LocalExistenceLemma}, let  $(u,v)$  be the local-in-time classical solution  of problem \eqref{problem}. If for some $\frac{n}{2}<p<n$ the $u$-component belongs  to $L^\infty((0,T_{max});L^p(\Omega))$, then $(u,v)$ is global in time, i.e. $T_{max}=\infty$, and moreover both $u$ and $v$ are bounded in $\Omega \times (0,\infty)$. 
\begin{proof}
For $u\in L^\infty((0,T_{max});L^p(\Omega))$, taking into consideration assumption \eqref{AssumptionsOng} on $g$, we have that 
\[
\int_\Omega g(u)^p\leq \lambda_2^p\int_\Omega (1+u)^{\beta p}\leq \lambda_2^p \int_\Omega (1+u)^p\quad \textrm{for all}\quad t\in(0,T_{max}),
\]
so that also $g(u)\in L^\infty((0,T_{max});L^p(\Omega))$. Consequently, standard elliptic regularity results applied to the second equation of \eqref{problem} warrant $v \in L^\infty((0,T_{max});W^{2,p}(\Omega))$ and hence $\nabla v\in  L^\infty((0,T_{max});W^{1,p}(\Omega))$, and finally Sobolev embedding theorems give $ v\in  L^\infty((0,T_{max});C^{[2-(n/p)]}(\bar{\Omega}))$ and $\nabla v\in  L^\infty((0,T_{max});L^{q}(\Omega))$ for all $n<q<p^*:=\frac{np}{n-p}$. In particular, for some positive constant $C_q$ we have that 
\begin{equation}\label{Bound_v_1-q}
\lVert  v (\cdot, t)\lVert_{L^{q}(\Omega)}+\lVert  \nabla v (\cdot, t)\lVert_{L^{q}(\Omega)}  \leq C_q\quad \textrm{for all}\quad t\in(0,T_{max}),
\end{equation}
and in addition, because $q>n$ implies $W^{1,q}(\Omega) \xhookrightarrow{} L^\infty(\Omega)$, we also attain $v\in L^\infty(\Omega \times (0,T_{max}))$.

As far as $u$ is concerned, for any $x\in \Omega$ and $t \in (0,T_{max})$ we set $t_0:=\max\{0,t-1\}$ so that  the representation formula for $u$ yields
\begin{equation}\label{RepresentationFormulaU}
\begin{split}
u (\cdot,t) &\leq e^{(t-t_0)\Delta}u(\cdot,t_0)-\int_{t_0}^t e^{(t-s)\Delta}\nabla \cdot (u (\cdot,s) \chi(\cdot,s)\nabla v (\cdot,s))ds\\ & 
=: u_{1}(\cdot,t)+u_{2}(\cdot,t). 
\end{split}
\end{equation}
Here, we invoke known smoothing estimates for the Neumann heat semigroup (see  \cite[Lemma 1.3]{WinklAggre}) which warrant the existence of positive constants $C_s$ and $\hat{C}_S$  such that for all $t>0$ and $f\in L^1(\Omega)$
 \begin{equation} \label{LpLqEstimate}   
 \lVert  e^{t\Delta} f  \lVert_{L^\infty(\Omega)}\leq C_S (1+t^{-\frac{n}{2}}) \lVert f\lVert_{L^1(\Omega)},
 \end{equation} 
 and for all $t>0$, $r>1$ and $f\,\in C^1(\bar{\Omega})$ with $\frac{\partial f}{\partial \nu}=0$ on $\partial \Omega$
 \begin{equation}\label{LpLqEstimateGradient}
 \lVert  e^{t\Delta}\nabla \cdot  \nabla f  \lVert_{L^\infty(\Omega)}\leq \hat{C}_S (1+t^{-\frac{1}{2}-\frac{n}{2r}})e^{-\lambda_1 t} \lVert \nabla f\lVert_{L^r(\Omega)},
 \end{equation}
 $\lambda_1>0$ denoting the first nonzero eigenvalue of $-\Delta$ in $\Omega$ under Neumann boundary conditions. 

Subsequently, if  $t\leq 1$ and hence $t_0=0$ we achieve from the parabolic maximum principle 
\begin{equation}\label{UnifrormBound_u1_t_smaller_one}
%\begin{cases}
\lVert u_{1}(\cdot,t)  \lVert_{L^\infty(\Omega)} \leq \lVert u_0 \lVert_{L^\infty(\Omega)}\quad \textrm{for all}\quad t\leq 1.
%\lVert u_{\varepsilon 1}(\cdot, )  \lVert_{L^\infty(\Omega)} \leq d_5\lVert u (\cdot,t_0) \lVert_{L^p(\Omega)}\leq d_5 C_p,&\;\;\textrm{for all}\;\,t> 1.\\
%\end{cases}
\end{equation}
Conversely, for all $t>1$ and hence $t-t_0=1$, an application of  \eqref{LpLqEstimate} with $f=u(\cdot,t_0)$ and bound \eqref{Bound_of_u} infer for all $1< t < T_{max}$ this estimate
\begin{equation}\label{UnifrormBound_u1_t_greather_one}
%\begin{split}
\lVert u_{ 1}(\cdot,t)  \lVert_{L^\infty(\Omega)} \leq C_S (1+(t-t_0)^\frac{-n}{2})\lVert u(\cdot,t_0) \lVert_{L^1(\Omega)}\leq 2 mC_S. %
%&\;\;\textrm{for all}\;\,t\leq 1,\\
%\lVert u_{\varepsilon 1}(\cdot, )  \lVert_{L^\infty(\Omega)} \leq d_5\lVert u (\cdot,t_0) \lVert_{L^p(\Omega)}\leq d_5 C_p,&\;\;\textrm{for all}\;\,t> 1.\\
%\end{split}
\end{equation}
Furthermore, for any $n<r<q$ we apply \eqref{LpLqEstimateGradient} with $\nabla f=u (\cdot,t) \chi(\cdot,t)\nabla v (\cdot,t)$ and arrive for $t\in (0,T_{max})$ at
\begin{equation}\label{UnifrormBound_u2_previous}
\begin{split}
&\lVert u_{2}(\cdot,t) \lVert_{L^\infty(\Omega)}  \leq  \int_{t_0}^t \lVert e^{(t-s)\Delta}\nabla \cdot (u (\cdot,s) \chi(\cdot,s)\nabla v (\cdot,s)) \lVert_{L^\infty(\Omega)}ds\\ &
\quad \leq \lVert \chi \rVert_{\Normchi}  \hat{C}_S \int_{t_0}^t (1+(t-s))^{-\frac{1}{2}-\frac{n}{2r}}\lVert u (\cdot,s) \nabla v (\cdot,s) \lVert_{L^r(\Omega)} ds,
%% \leq \\ & 
%%\quad c \Gamma \big(\frac{1}{2}-\zeta\big),
\end{split} 
\end{equation}
where once we invoked the boundedness of $\chi$ given in \eqref{MaximumPricnipleRelation}. 

Now, for any given $t'\in (0,T_{max})$ we consider the function defined by
\begin{equation}\label{DefinitionSequenceAEpsilon}
A(t'):=\sup_{t\in (0,t')}\lVert u (\cdot, t)\lVert_{L^\infty(\Omega)},
\end{equation}
which is bounded in view of the properties of $u$. Hence, an interpolation inequality, the assumption $u\in L^\infty((0,T_{max});L^p(\Omega))$ (i.e. $\int_\Omega u^p\leq C_p$ for some $C_p>0$ and all $t\in(0,T_{max})$) and \eqref{Bound_v_1-q} entail for all $s\,\in (0,t')$
\begin{equation*}
\begin{split}
\lVert u (\cdot,s) \nabla v (\cdot,s) \lVert_{L^r(\Omega)} &\leq \lVert u (\cdot,s) \lVert_{L^{\frac{rq}{q-r}}(\Omega)} \lVert \nabla v (\cdot,s) \lVert_{L^{q}(\Omega)} \\ &
\leq \lVert u (\cdot,s) \lVert_{L^{\infty}(\Omega)}^{1-\frac{p(q-r)}{rq}} \lVert u (\cdot,s) \lVert_{L^{p}(\Omega)}^{\frac{p(q-r)}{rq}} \lVert \nabla v (\cdot,s) \lVert_{L^{q}(\Omega)}\\ &
\leq A^{1-\frac{p(q-r)}{rq}} (t')C_p^{\frac{q-r}{rq}}C_q=c_4A^{l} (t'), 
\end{split}
\end{equation*}
with $c_4=C_p^{\frac{q-r}{rq}}C_q$ and $0<l=1-\frac{p(q-r)}{rq} <1$. Subsequently, in light of the estimate now gained for $u (\cdot,s) \nabla v (\cdot,s)$ and the relation $t-t_0\leq 1$, bound \eqref{UnifrormBound_u2_previous} reads
\begin{equation}\label{UnifrormBound_u2}
\begin{split}
\lVert u_{2}(\cdot, t) \lVert_{L^\infty(\Omega)} & \leq  \frac{2r}{r-n}  \lVert \chi \rVert_{\Normchi}  \hat{C}_S c_4 (2^{\frac{1}{2}-\frac{n}{2r}}-1)A^{l} (t')=:c_5A^{l}(t').
%\chi \int_{t_0}^t \lVert e^{(t-s)\Delta}\nabla \cdot (u (\cdot,s) \nabla v (\cdot,s)) \lVert_{L^\infty(\Omega)}ds\\ &
%\leq \chi C_Sd_8\int_{t_0}^t (1+(t-s))^{-\frac{1}{2}-\frac{n}{2r}}e^{-\mu_1 t}\lVertu (\cdot,s) \nabla v (\cdot,s) \lVert_{L^r(\Omega)} ds.
% \leq \\ & 
%\quad c \Gamma \big(\frac{1}{2}-\zeta\big),
\end{split} 
\end{equation}
From expression \eqref{RepresentationFormulaU}, by collecting  \eqref{UnifrormBound_u1_t_smaller_one}-\eqref{UnifrormBound_u2_previous} and \eqref{UnifrormBound_u2} we infer
\begin{equation*}
\begin{split}
\lVert u (\cdot,t) \lVert_{L^\infty(\Omega)}\leq c_6(A^{l}(t')+1)\quad \textrm{for all}\quad t \in \, (0,t'),
%\lVert e^{(t-t_0)\Delta}u(\cdot,t_0)\lVert_{L^\infty(\Omega)}\\ & \quad+\chi\int_{t_0}^t \lVert e^{(t-s)\Delta}\nabla \cdot (u (\cdot,s) \nabla v (\cdot,s))\lVert_{L^\infty(\Omega)}ds \\ &
%\quad +\int_{t_0}^t \lVert e^{(t-s)\Delta}a\lVert_{L^\infty(\Omega)}ds.
\end{split}
\end{equation*}
where $c_6=\max\{\max\{||u_0||_{L^\infty(\Omega)},2mC_S\},c_5\}$. Therefore recalling \eqref{DefinitionSequenceAEpsilon} 
\begin{equation*}
\sup_{t\in (0,t')}\lVert u (\cdot, t)\lVert_{L^\infty(\Omega)}=:A(t')\leq c_6(A^{l}(t')+1)\quad \textrm{for all}\quad t'\in (0,T_{max}),
\end{equation*}
which through Lemma \ref{BoundsInequality} yields this bound for $u$:
\begin{equation*} %\label{BoundUnfiromuepsilon}
\lVert u (\cdot,t) \lVert_{L^\infty(\Omega)}\leq \max\{1,(2c_6)^\frac{1}{1-l}\}=:\hat{c} \quad \textrm{for all}\quad t\in (0,t').
\end{equation*}
Since $\hat{c}>0$ is time-independent and $t'\in(0,T_{max})$ is arbitrary, the above uniform bound for $u$ holds up $T_{max}$. Hence the extensibility criterion \eqref{extensibility_criterion_Eq} of Lemma \ref{LocalExistenceLemma} shows that $T_{max} = \infty$ and that both $u$ and $v$ are bounded in $\Omega \times (0,\infty)$. 
\end{proof}
\end{lemma}
\section{A priori estimates and proof of the main result}\label{EstimatesAndProofSection}
In this section we shall gain some uniform bound for $u$, by deriving an upper bound for $\lVert u\lVert_{L^p(\Omega)}$, with $p$ sufficiently large and on the whole interval $(0,T_{max})$, which is given by a positive and time independent constant. This is attained by constructing an absorptive differential inequality for $t\mapsto \int_\Omega u^p$ and using comparison principles, exactly as specified in this sequel of lemmas.
\begin{lemma}\label{Estim_GradientvSquareLemma} 
Under the assumptions of Lemma \ref{LocalExistenceLemma}, let $(u,v)$  be the local-in-time classical solution of problem \eqref{problem}. Then there exists a positive constant $M$ such that 
\begin{equation}\label{Estim_GradientVSquare}   
\begin{split}
\int_\Omega|\nabla v|^2\leq M\quad \textrm{for all}\quad t\in (0,T_{max}). %\\ &
%\quad+ p(p-1)\chi_0 C_1(\epsilon_1)\epsilon_2 2 (4p^2+n)\lVert v_0vrt^2_{L^\infty(\Omega)}\int_\Omega \lvert \nabla  v\rvert^{2p-2} \lvert D^2 v\rvert^2\\ &
%\quad +|\Omega|p((2\mu+k_+)C_3(\epsilon_3)+c_0(p-1)\chi_0),
\end{split}
\end{equation}
\begin{proof}
Testing the second equation of  \eqref{problem} by $v$, we obtain
\[
-\int_\Omega v \Delta v=-\int_\Omega v^2+\int_\Omega g(u)v \quad \textrm{for all}\quad t\in (0,T_{max}),
\]
so that integrations by part, the mass conservation property \eqref{Bound_of_u} and the Young inequality (recall also \eqref{AssumptionsOng}) lead for all $t \in (0, T_{max})$ to
\[
\begin{split}
\int_\Omega |\nabla v|^2& \leq -\int_\Omega v^2+\int_\Omega v^2+\frac{\lambda_2^2}{4}\int_\Omega (1+u)^{2\beta}
\leq \frac{\lambda_2^2}{4}(|\Omega|+m)=:M.
%\begin{cases}
%\frac{\lambda_2^2}{4}(2\beta (|\Omega|+m)+(1-2\beta)|\Omega|)& \textrm{if}\quad 0\leq\beta<\frac{1}{2},  \\
%\frac{\lambda_2^2}{4}(|\Omega|+m)  &\textrm{if}\quad \beta=\frac{1}{2}.
%\end{cases}
\end{split}
\]
\end{proof}
\end{lemma}
\begin{lemma}\label{Estim_general_For_u^pLemma}  
Under the assumptions of Lemma \ref{LocalExistenceLemma}, let $(u,v)$ be the local-in-time classical solution  of problem \eqref{problem}. Then, for any $p>1$
\begin{equation}\label{Estim_general_For_u^p}   
\begin{split}
\frac{d}{dt}\int_\Omega u^p&+p\frac{p-1}{2}\int_\Omega u^{p-2}\lvert \nabla u\rvert^2 \leq \\ &
+p\frac{p-1}{2}\lVert \chi \rVert^2_{L^\infty([\gamma,\Gamma])}\int_\Omega u^{p}|\nabla v|^2\quad \textrm{on}\quad (0,T_{max}). %\\ &
%\quad+ p(p-1)\chi_0 C_1(\epsilon_1)\epsilon_2 2 (4p^2+n)\lVert v_0vrt^2_{L^\infty(\Omega)}\int_\Omega \lvert \nabla  v\rvert^{2p-2} \lvert D^2 v\rvert^2\\ &
%\quad +|\Omega|p((2\mu+k_+)C_3(\epsilon_3)+c_0(p-1)\chi_0),
\end{split}
\end{equation}
%where $C_1=\frac{\lambda_2}{p}|\Omega|(\lambda_2(p-1))^{p-1}$.
%\[
%\begin{cases}
%%C_1(\epsilon_1)=\frac{1}{4 \epsilon_1}\quad C_2(\epsilon_2)=\frac{p}{p+1} (\epsilon_2(p+1))^\frac{-1}{p} \quad    
% C_1(\epsilon)=\frac{1-\alpha}{p+1}\left(\frac{\epsilon(p+1)}{2p(p+\alpha)}\frac{p+\alpha-1}{p -1}\right)^\frac{p+\alpha}{\alpha-1},
% \\C_2(\epsilon)=\frac{1}{p+1}(\frac{\epsilon(p+1)\chi_0}{2(2\mu+k) p^2})^{-p},
% \\ c_0=\frac{p-1}{p+\alpha-1}\chi_0C_1(\epsilon)|\Omega|+(2\mu+k)C_2(\epsilon)|\Omega|.
%\end{cases}
%\]
%$\hat{C}_{GN}=\chi^2p\frac{p-1}{2(p+1)}C_{GN}$, with $C_{GN}$ as in Lemma \ref{InequalityG-NLemma} and $C(\sigma_1)$ a $\sigma_1$-depending constant, $\sigma_1$  being a positive constant to be fixed later on.
\begin{proof}
For $p>1$, testing the first equation of problem \eqref{problem} by $u^{p-1}$ and using its boundary conditions  provide
\begin{equation}\label{Estim_1_For_u^p}  
\begin{split}
\frac{1}{p}\frac{d}{dt} \int_\Omega u^p&=\int_\Omega u^{p-1}u_t =-(p-1)\int_\Omega u^{p-2}\lvert \nabla u\rvert^2  \\ &\quad +(p-1) \int_\Omega  u^{p-1}\chi(v)\nabla u \cdot  \nabla v \quad \textrm{for all}\quad t\in (0,T_{max}).
% \\ &\quad+\lambda_2\int_\Omega (u+1)^{p-1}-\mu_2 \int_\Omega u(u+1)^{p-1}, 
\end{split}
\end{equation}
As to the integral involving $u^{p-1}\chi(v)\nabla u \cdot  \nabla v$, the Young inequality and the uniform bound for $\chi$ derived in the third relation of \eqref{MaximumPricnipleRelation}  yield for all $t\in(0,T_{max})$
\[
\begin{split}
\int_\Omega|u^{p-1}\chi(v)\nabla u \cdot  \nabla v|& \leq \frac{1}{2} \int_\Omega u^{p-2}\chi^2(v)|\nabla u|^2+ \frac{1}{2}\int_\Omega u^{p}|\nabla v|^2 \\ &
\leq  \frac{1}{2}\lVert \chi\rVert_{\Normchi}^2\int_\Omega u^{p-2}|\nabla u|^2+\frac{1}{2}\int_\Omega u^{p}|\nabla v|^2,
\end{split}
\]
so that equality \eqref{Estim_1_For_u^p}  writes exactly as claimed.
\end{proof}
\end{lemma}
\begin{lemma}\label{Estim_general_For_TermWithChiSquareLemma} 
For  $n=2$ and under the remaining assumptions of Lemma \ref{LocalExistenceLemma}, let $(u,v)$ be the local-in-time classical solution  of problem \eqref{problem}. Then for any $p\in(1,2)$ and all  $t\in (0,T_{max})$ holds 
\begin{equation}\label{MainExpressionForFirstDerivativeuPowerP}
\frac{d}{d t}\int_\Omega u^p \leq -\frac{p(p-1)}{2}\int_\Omega u^{p-2}|\nabla u|^2+(\varepsilon+c_8\delta)\int_\Omega u^{p+1}+c_{10}, 
\end{equation}
where $c_8$ and $c_{10}$ are computable positive constants, $\varepsilon$ an arbitrary  positive number and $\delta$ another arbitrary nonnegative real such that 
\[
\begin{cases}
\delta=0 & \textrm{if}\quad \beta=0,\\
\delta> 0 & \textrm{if}\quad \beta>0.\\
\end{cases}
\] 
\begin{proof}
In order to estimate the term $\int_\Omega u^p |\nabla v|^2$ appearing in the above Lemma \ref{Estim_general_For_u^pLemma}, we observe that an application of the Hölder inequality infers for any $p>1$ 
\begin{equation}\label{EstimateFirstu^pnablav^2}
\int_\Omega u^p |\nabla v|^2 \leq \bigg(\int_\Omega u^{p+1}\bigg)^\frac{p}{p+1}\bigg(\int_\Omega |\nabla v|^{2p+2}\bigg)^\frac{1}{p+1}=\lVert u \rVert^p_{L^{p+1}(\Omega)}\lVert \nabla v \rVert^2_{L^{2p+2}(\Omega)}.
\end{equation} 
Now, in view of the fact that the $v$-component solves $-\Delta v+ v=g(u)$ in $\Omega\times (0,T_{max})$,  for $p\in(1,2)$ we can invoke Lemma \ref{LemmaTwoDimensional} with $f=v$ which, together with \eqref{Estim_GradientVSquare},  implies
\begin{equation}\label{EstimateSecondu^pnablav^2}
\begin{split}
\bigg(\int_\Omega |\nabla v|^{2p+2}\bigg)^\frac{1}{p+1}&=\lVert \nabla v \rVert^2_{L^{2p+2}(\Omega)}\leq \hat{C}^2 \lVert -\Delta v + v\rVert_{L^{p+1}(\Omega)} \lVert \nabla v\rVert_{L^2(\Omega)}\\ &
\leq \hat{C}^2\sqrt{M}  \lVert g(u)\rVert_{L^{p+1}(\Omega)}\quad \textrm{for all}\quad t\in (0,T_{max}).
\end{split}
\end{equation}
Moreover, assumption \eqref{AssumptionsOng} gives
\begin{equation}\label{EstimateThirdu^pnablav^2}
\lVert g(u)\rVert_{L^{p+1}(\Omega)}\leq \lambda_2 \bigg(\int_\Omega (1+u)^{\beta(p+1)}\bigg)^\frac{1}{p+1} \quad \textrm{on}\quad (0,T_{max}),
\end{equation}
so that, by means of Young's inequality and the introduction of a positive real number $\varepsilon$, we deduce from \eqref{EstimateFirstu^pnablav^2}, \eqref{EstimateSecondu^pnablav^2},  \eqref{EstimateThirdu^pnablav^2} and \eqref{AlgebraicInequality2toalpha} that on  $(0,T_{max})$ 
\begin{equation}\label{EstimateFouru^pnablav^2}
\begin{split}
\int_\Omega u^p |\nabla v|^2 &\leq \bigg(\int_\Omega u^{p+1}\bigg)^\frac{p}{p+1}\bigg((\hat{C}^2\sqrt{M}\lambda_2)^{p+1}\int_\Omega (1+u)^{\beta(p+1)}\bigg)^\frac{1}{p+1}\\ &
\leq \frac{\varepsilon}{c_8}\int_\Omega u^{p+1} +\frac{c_7}{p+1}\bigg(\frac{\varepsilon(p+1)}{c_8p}\bigg)^{-p}\int_\Omega (1+u)^{\beta(p+1)}\\ &
\leq 
\begin{cases}
\frac{\varepsilon}{c_8}\int_\Omega u^{p+1} +c_92^{\beta(p+1)}|\Omega|+c_92^{\beta(p+1)}\int_\Omega u^{\beta(p+1)},& \beta >0,\\
\frac{\varepsilon}{c_8}\int_\Omega u^{p+1} +c_9|\Omega|& \beta =0,
\end{cases}
\end{split}
\end{equation} 
with $c_7=(\hat{C}^2\sqrt{M}\lambda_2)^{p+1}$, $c_8=\frac{p(p-1)}{2}\lVert \chi \rVert_{\Normchi}^2$ and $c_9=\frac{c_7}{p+1}\big(\frac{\varepsilon(p+1)}{c_8p}\big)^{-p}$. In particular, for $\beta >0$, similar fashions allow us to show that for some $\delta>0$, on $(0,T_{max})$ also holds 
\begin{equation}\label{EstimateFiveu^pnablav^2}
\int_\Omega   u^{\beta(p+1)}\leq \frac{\delta}{c_9 2^{\beta(p+1)}}\int_\Omega u^{p+1} + (1-\beta)\bigg(\frac{\delta}{c_9 \beta 2^{\beta(p+1)}}\bigg)^{-\frac{\beta}{1-\beta}}|\Omega|.
\end{equation}
Finally, taking into account \eqref{Estim_general_For_u^p}  of Lemma \ref{Estim_general_For_u^pLemma}, with relations  \eqref{EstimateFouru^pnablav^2} and \eqref{EstimateFiveu^pnablav^2}  we conclude the proof setting 
\[
c_{10}=
\begin{cases}
c_8c_9 |\Omega| & \textrm{for}\quad \beta=0,\\
c_8c_9 2^{\beta(p+1)}|\Omega|(1+(1-\beta)\big(\frac{\delta}{c_9 \beta 2^{\beta(p+1)}}\big)^{-\frac{\beta}{1-\beta}})&\textrm{for}\quad \beta>0.
\end{cases}
\]
\end{proof}
\end{lemma}
As announced, the succeeding lemma is dedicated to establish the desired uniform-in-time bound for $\lVert u\lVert_{L^p(\Omega)}$, with some proper $p>1$. 
\begin{lemma}\label{LemmaAbsorptiveMainInequality}
For  $n=2$ and under the remaining assumptions of Lemma \ref{LocalExistenceLemma}, let $(u,v)$ be the local-in-time classical solution  of problem \eqref{problem}. Then for any $p\in(1,2)$ there exists a positive constant $C_p$ such that 
\begin{equation}\label{BoundU+1EnadNabla}
\int_\Omega u^p \leq C_p \quad \textrm{for all}\quad t \in (0,T_{max}).
\end{equation}
\begin{proof}
We rely on the Gagliardo-Nirenberg inequality to estimate the term $\int_\Omega u^{p+1}$ appearing in \eqref{MainExpressionForFirstDerivativeuPowerP}. Precisely, expression \eqref{InequalityTipoG-N} with $f=u^\frac{p}{2}$, $\mathfrak{j}=0, \mathfrak{m}=1, \mathfrak{r}=n=2, \mathfrak{p}=\frac{2(p+1)}{p}$ and $\mathfrak{q}=\mathfrak{s}=\frac{2}{p}$, in conjunction with \eqref{AlgebraicInequality2toalpha}, provides  for  all $t\in (0,T_{max})$ the relation
\begin{equation*} %\label{Estim_3_For_u^p}
\begin{split}
\int_\Omega u^{p+1}&=\lvert \lvert u^\frac{p}{2}\lvert \lvert_{L^\frac{2(p+1)}{p}(\Omega)}^\frac{2(p+1)}{p}\leq c_{11}        \lvert \lvert\nabla u^\frac{p}{2}\lvert \lvert_{L^2(\Omega)}^{\frac{2(p+1)}{p}\theta_1}      \lvert \lvert u^\frac{p}{2}\lvert \lvert_{L^\frac{2}{p}(\Omega)}^{\frac{2(p+1)}{p}(1-\theta_1)} +c_{11}      \lvert \lvert u^\frac{p}{2}\lvert \lvert^\frac{2(p+1)}{p}_{L^\frac{2}{p}(\Omega)},
 \end{split}
\end{equation*}
being $c_{11}= (2C_{GN})^\frac{2(p+1)}{p}$ and $0<\theta_1=\frac{p}{p+1}<1$. In particular, this gained inequality and \eqref{Bound_of_u} give
\begin{equation*} %\label{Estim_3_For_u^p}
\begin{split}
(\varepsilon+c_8\delta)\int_\Omega u^{p+1}&\leq(\varepsilon+c_8\delta) c_{11} m \int_\Omega |\nabla u^\frac{p}{2}|^2 \\ & 
\quad + (\varepsilon+c_8\delta)c_{11}m^{p+1} \quad \textrm{for all}\quad t \in (0,T_{max}),
 \end{split}
\end{equation*}
so that \eqref{MainExpressionForFirstDerivativeuPowerP} is transformed in 
 \begin{equation}\label{Derivativeu^pSecondstep}
\frac{d}{d t}\int_\Omega u^p \leq \Bigg(c_{11}m(\varepsilon+c_8\delta)-\frac{2(p-1)}{p}\Bigg)\int_\Omega |\nabla u^\frac{p}{2}|^2 + c_{12}, 
\end{equation}
with  $c_{12}=(\varepsilon+c_8\delta)c_{11}m^{p+1}+c_{10}$. 

Successively we again use Lemma \ref{InequalityG-NLemma}  and in the Gagliardo-Nirenberg inequality  \eqref{InequalityTipoG-N} we take $f=u^\frac{p}{2}$, $\mathfrak{q}=\mathfrak{s}=\frac{2}{p}$, $\mathfrak{p}=n=\mathfrak{r}=2$, $\mathfrak{j}=0$ and $\mathfrak{m}=1$; we infer through \eqref{AlgebraicInequality2toalpha} that on $(0,T_{max})$ one has
\begin{equation}\label{Estim_3_For_u^p}
\begin{split}
\int_\Omega u^{p}&=\lvert \lvert u^\frac{p}{2}\lvert \lvert_{L^2(\Omega)}^2\leq c_{13}       \lvert \lvert\nabla u^\frac{p}{2}\lvert \lvert_{L^2(\Omega)}^{2\theta_2}   \lvert \lvert u^\frac{p}{2}\lvert \lvert_{L^\frac{2}{p}(\Omega)}^{2(1-\theta_2)} +c_{13}   \lvert \lvert u^\frac{p}{2}\lvert \lvert^2_{L^\frac{2}{p}(\Omega)},
 \end{split}
\end{equation}
where $0<\theta_2=\frac{p-1}{p}<1$ and $c_{13}= (2C_{GN})^2$. Considering once again bound \eqref{Bound_of_u} and  introducing $c_{14}=c_{13} \max\{m,m^p\}$, by making again use of \eqref{AlgebraicInequality2toalpha}  inequality \eqref{Estim_3_For_u^p} can also be rewritten as 
\begin{equation}\label{EstimGradientu_pOver2}
\begin{split}
-\int_\Omega \lvert \nabla u^\frac{p}{2}\rvert^2\leq 1-(2c_{14})^\frac{p}{1-p}\left(\int_\Omega u^{p}\right)^\frac{p}{p-1} \quad \textrm{for all} \quad t \in (0,T_{max}).
 \end{split} 
\end{equation}
We, then, choose 
\[
\begin{cases}
\textrm{if}\quad \beta=0, \quad (\textrm{and so}\quad \delta=0)& \varepsilon=\frac{p-1}{c_{11}mp}>0,\\
\textrm{if}\quad \beta >0, & \varepsilon=\frac{p-1}{2c_{11}mp}>0 \;\;\textrm{and}\;\; \delta=\frac{p-1}{2c_8c_{11}mp}>0,
\end{cases}
\] 
and plug such values and  estimate \eqref{EstimGradientu_pOver2} into \eqref{Derivativeu^pSecondstep}; in this way we arrive at this initial problem
\begin{equation*}\label{MainInitialProblemWithM}
\begin{cases}
\Phi'(t)\leq c_{15}-c_{16} \Phi(t)^\frac{p}{p -1}\quad t \in (0,T_{max}),\\
\Phi(0)=\int_\Omega u_0^p, 
\end{cases}
\end{equation*}
where $\Phi(t):=\int_\Omega u^p$, $c_{15}=c_{12}+\frac{p-1}{p}$ and $c_{16}=\frac{p-1}{p}(2c_{14})^\frac{p}{1-p}$. Ultimately, we conclude the demonstration by an application of a comparison principle implying 
\begin{equation*} %\label{BoundU+1^pFprK=1}
\Phi(t)\leq \max\left\{\Phi(0),\left(\frac{c_{15}}{c_{16}}\right)^\frac{p -1}{p}\right\}=:C_p\quad \textrm{for all}\quad t\in(0,T_{max}).
\end{equation*}
\end{proof}
\end{lemma}
As a consequence of all of the above preparations, we finally can prove our claimed statement:
\\ \\
\textit{Proof of Theorem \ref{MainTheorem}}. For $n=2$  and any nonnegative $0\not \equiv u_0\in C^0(\bar{\Omega})$, Lemma \ref{LocalExistenceLemma} provides a unique local-in-time classical solution $(u,v)$ to problem \eqref{problem}. Thereafter, by virtue of Lemma \ref{LemmaAbsorptiveMainInequality}, for all $p\in(1,2)$ relation \eqref{BoundU+1EnadNabla} is achieved, then $u\in L^\infty((0,T_{max}); L^p(\Omega))$ and we can conclude through Lemma \ref{FromLocalToGLobalBoundedLemma}.
\qed
% % % % %
\section{Numerical simulations}\label{SimuationsSection}
In this section we numerically test the presented theoretical results by simulating system \eqref{problem} in two dimensions; for simplicity and with no possibility of confusion, the spatial variable $x=(x_1,x_2)$ is indicated with $(x,y)$. Further, we investigate whether the solutions are: stationary (time independent) and homogeneous (spatially uniform density); stationary and heterogeneous (spatially non-uniform density), but bounded; or suffer from chemotactic blow-up in finite time.

Specifically, we use finite element methods to simulate system \eqref{problem} with a variety of functions defining $\chi(v)$ and $g(u)$. The solution algorithm is based on an adaptive, implicit Runge-Kutta finite element method (see \cite{Ascher-1997-IER}). The space $\Omega$ is defined to be the interior of the $[0,0.1]\times[0,0.1]$ square, making $\partial\Omega$ the boundary of the square. To some extent the size of the space is arbitrary since we can use a larger domain, rescale the system with respect to the size and, thus, reproduce the same dynamics. However, because of anticipating large and thin spike structures in the cells' density we choose to simulate a small domain, which allows us to resolve such small heterogeneous solutions more accurately, without increasing the overall simulation mesh resolution.

From system \eqref{problem} we impose to the unknowns $u$ and $v$ to satisfy Neumann conditions on the boundary. The initial condition for $u$ was chosen to be uniform constant, $\bar{u}$, plus noise, that is 
\begin{equation}
u(x,y,0)=u_0=|\bar{u}+\sigma\eta(x,y)\nonumber|,
\end{equation}
where $\eta:\Omega\rightarrow [-1/2,1/2]$ is a continuous uniformly random variable and $\sigma$ is a positive constant that allows the stochastic perturbation to be scaled with respect to the mean value $\bar{u}$. On the other hand, even though $v(x,y,0)=v_0$ is not required in system \eqref{problem} (exactly because the equation for $v$ is elliptic), in order to run the Runge-Kutta iterative method an initial condition for $v$ has to be also assigned: in particular, we take $v_0$ as the solution of $-\Delta v_0+v_0=g(u_0)$ under Neumann boundary conditions. 

Since we are looking for differences between solutions that are bounded, versus those that suffer from blow-up, if a solution is found to increase indefinitely, then the same simulation was repeated with a finer discretisation to ensure that this outcome is the true numerical solution, rather than a numerical artefact. Specifically, whenever a solution was observed to be grow without bound, the grid was refined to have ten times as many elements as previously simulated, to ensure the outcome. As we will see later, system \eqref{problem} can support heterogeneous spike solutions in the densities of $u$. Critically, as parameters of interest are altered, the spikes densities become larger, whilst their support becomes smaller. As the spikes tend to the form of a $\delta$-function, numerical solution convergence requires finer discretisations in order to resolve the spike's shape. However, limitations in computer memory and processing power limit this procedure of refinement. Such cases will be clearly highlighted.

\begin{figure}[h!!!t!!!b!!!]
\centering
\includegraphics[width=\textwidth]{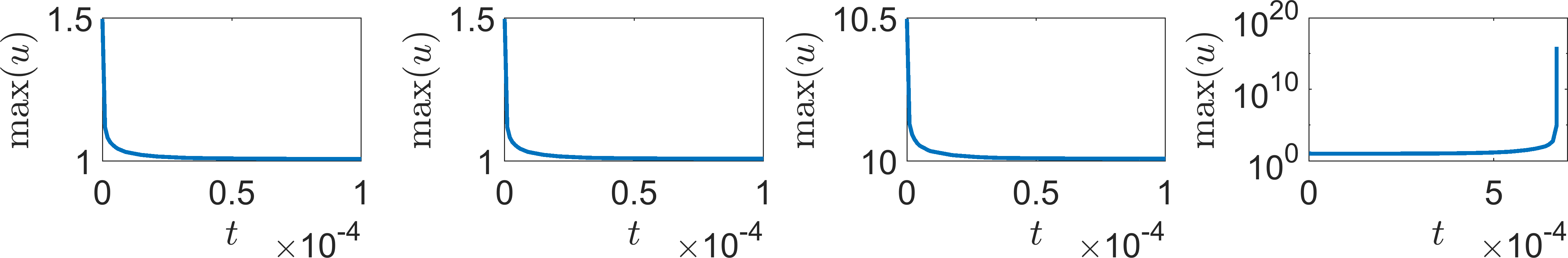}
\caption{\label{Base_case} Tracking the maximum value of $u$ on $\bar{\Omega}$ from simulations of system \eqref{problem} with $g(u)=u$ and $\chi(v)=\chi$; it is seen that the data  $\bar{u}$ and $\chi$ have to produce a large enough value of $\chi \int_\Omega u_0$ for blow-up to occur. From the left, the initial conditions are $\bar{u}=1, 1, 10, 10$ and $\chi=10^2, 10^3, 10^2, 10^3$. Throughout all simulations $\sigma=1$. Note that the scales on the first three plots are uniformly scaled, while the fourth plot (the right most) uses a logarithmic scale for the vertical axis (i.e. $\max (u)$), exactly to illustrate the unbounded growth of the solution.}
\end{figure}

Figure \ref{Base_case} illustrates one of the results discussed in $\S$\ref{Intro}. To be precise, when $\chi(v)=\chi>0$, constant, and $g(u)=u$ then solutions blow up if the initial mass $\int_\Omega u_0$ and $\chi$ are such that $\chi \int_\Omega u_0$ is  sufficiently large. Critically, when $\chi=10^3$ and $\bar{u}=10$ (right image of Figure \ref{Base_case}), then $u$ grows to over $10^{10}$ in less than $10^{-5}$ time units. The three simulations from the left of Figure \ref{Base_case} show what happens when either (or both) the initial condition, $\bar{u}$, or the sensitivity coefficient, $\chi$, is reduced by a factor of 10, namely the solution rapidly converges to the homogeneous steady state,
\begin{equation}
(u_s,v_s)=\left( \frac{1}{|\Omega|}\int_\Omega u(x,y,0),g(u_s)\right).\nonumber
\end{equation}

\begin{figure}[h!!!t!!!b!!!]
\centering
\subfigure[]{\includegraphics[width=.32\textwidth]{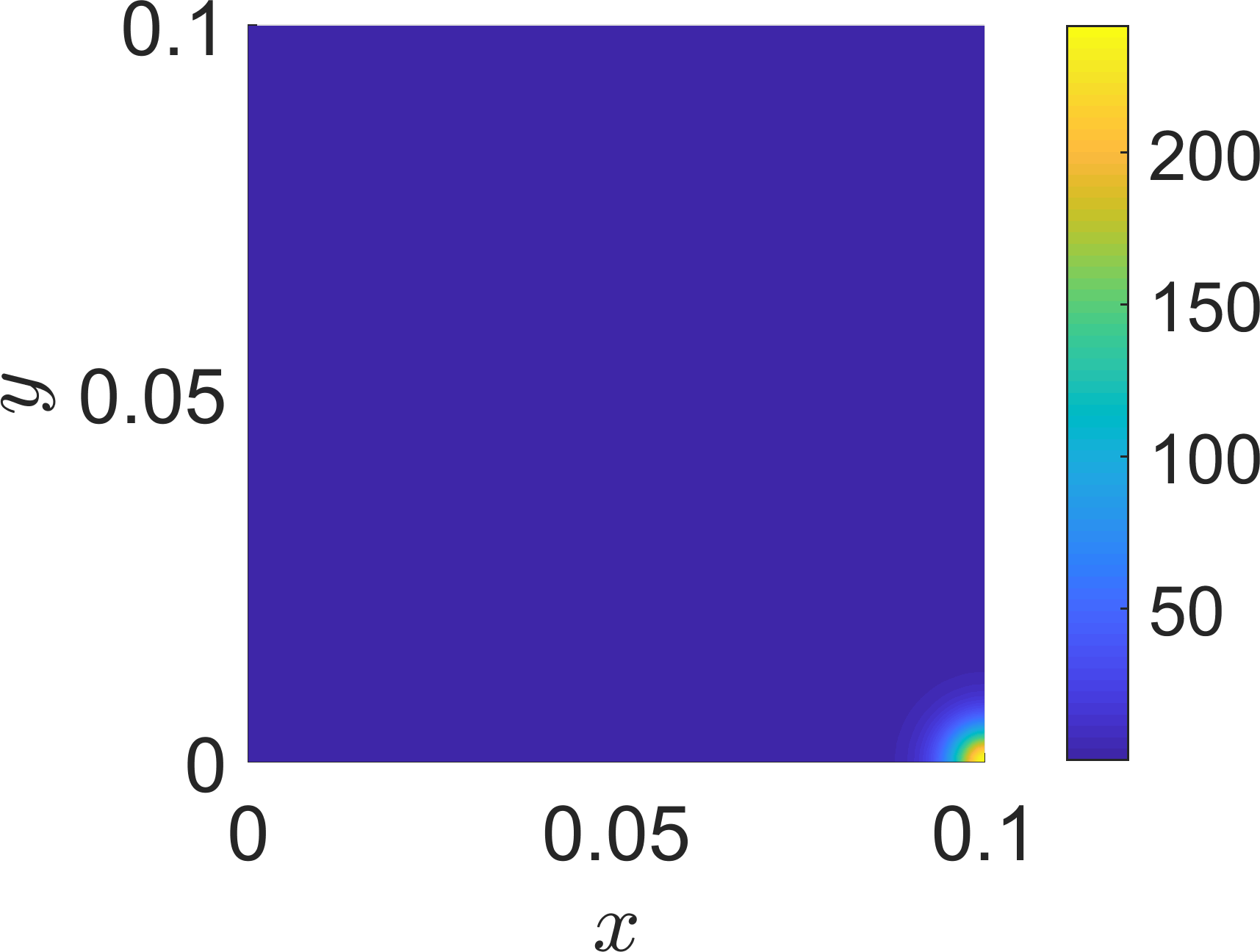}}
\subfigure[]{\includegraphics[width=.32\textwidth]{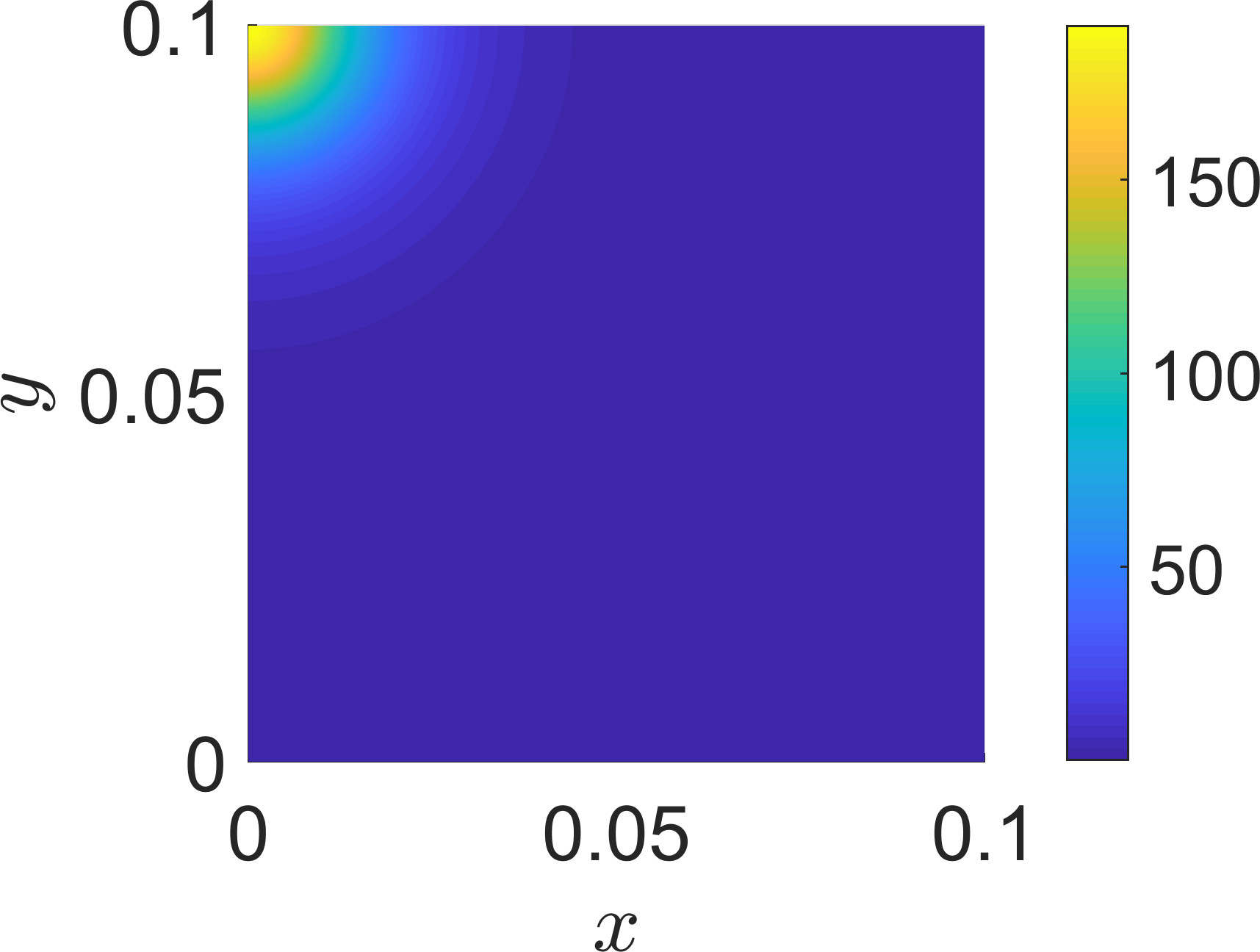}}
\subfigure[]{\includegraphics[width=.32\textwidth]{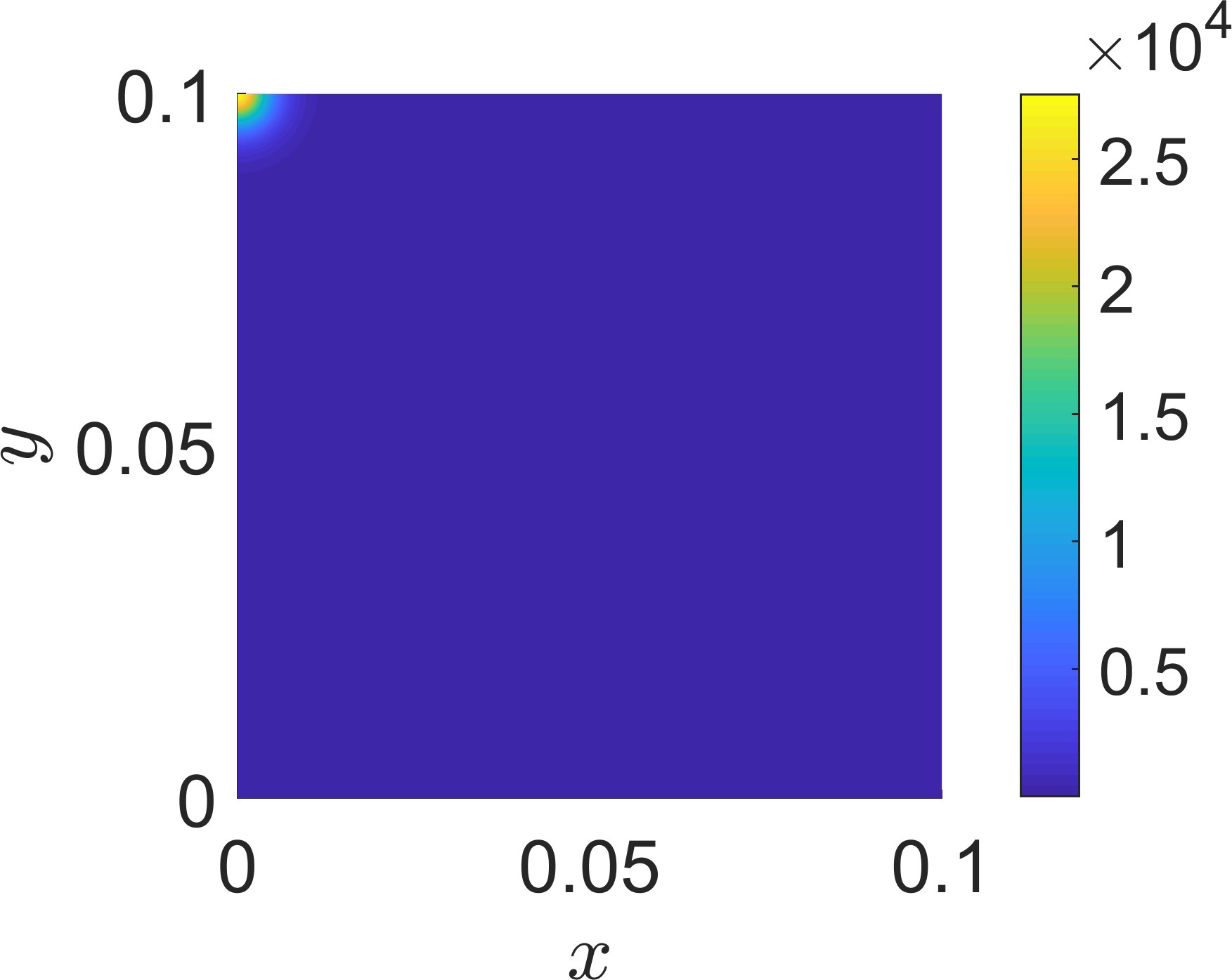}}
\caption{\label{g_sqrtu} Simulating system \eqref{problem} with $g(u)=(1+u)^{1/2}$ and $\chi(v)=\chi$ and indicating that blow-up does not appear to occur, but heterogeneous spatial solutions are possible. Parameters are (a) $\bar{u}=\sigma=1$, $\chi=10^4$, (b) $\bar{u}=\sigma=10$, $\chi=10^3$ and (c) $\bar{u}=\sigma=100$, $\chi=10^3$. Due to differing scales each subfigure has its own colour bar specifying the solution density. Each image is taken at $t=100$, at which point the simulations are seen to have stopped evolving.}
\end{figure}

In Figure \ref{g_sqrtu} we simulate system \eqref{problem} again, but this time with $g(u)=(1+u)^{1/2}$. Initially, we simulated this new system using the same parameters, $\bar{u}$, $\sigma$ and $\chi$, as those specified in Figure \ref{Base_case} (data not shown) and quickly realised that the simulations no longer succumbed to blow-up, rather each variable $(u,v)$ converged to a finite stationary distribution. Figure \ref{g_sqrtu} demonstrates we are able to increase $\bar{u}$ and $\chi$ without fear of blow-up. However, increasing these two parameters leads to the uniform steady stable being driven unstable, thus, the system evolves to a stable, bounded heterogeneous density.

 As specified above, increasing the initial condition, $\bar{u}$, and/or the value of the sensitivity, $\chi$, causes the spike solution to become sharper, giving problems with numerical convergence as the values are increased. Thus, although the simulated solutions begin to grow as $\bar{u}$ and $\chi$ are enlarged, this is probably an issue of the numerical resolution, rather than an analytical singularity.

\begin{figure}[h!!!t!!!b!!!]
\centering
\subfigure[\label{homo}]{\includegraphics[width=.4\textwidth]{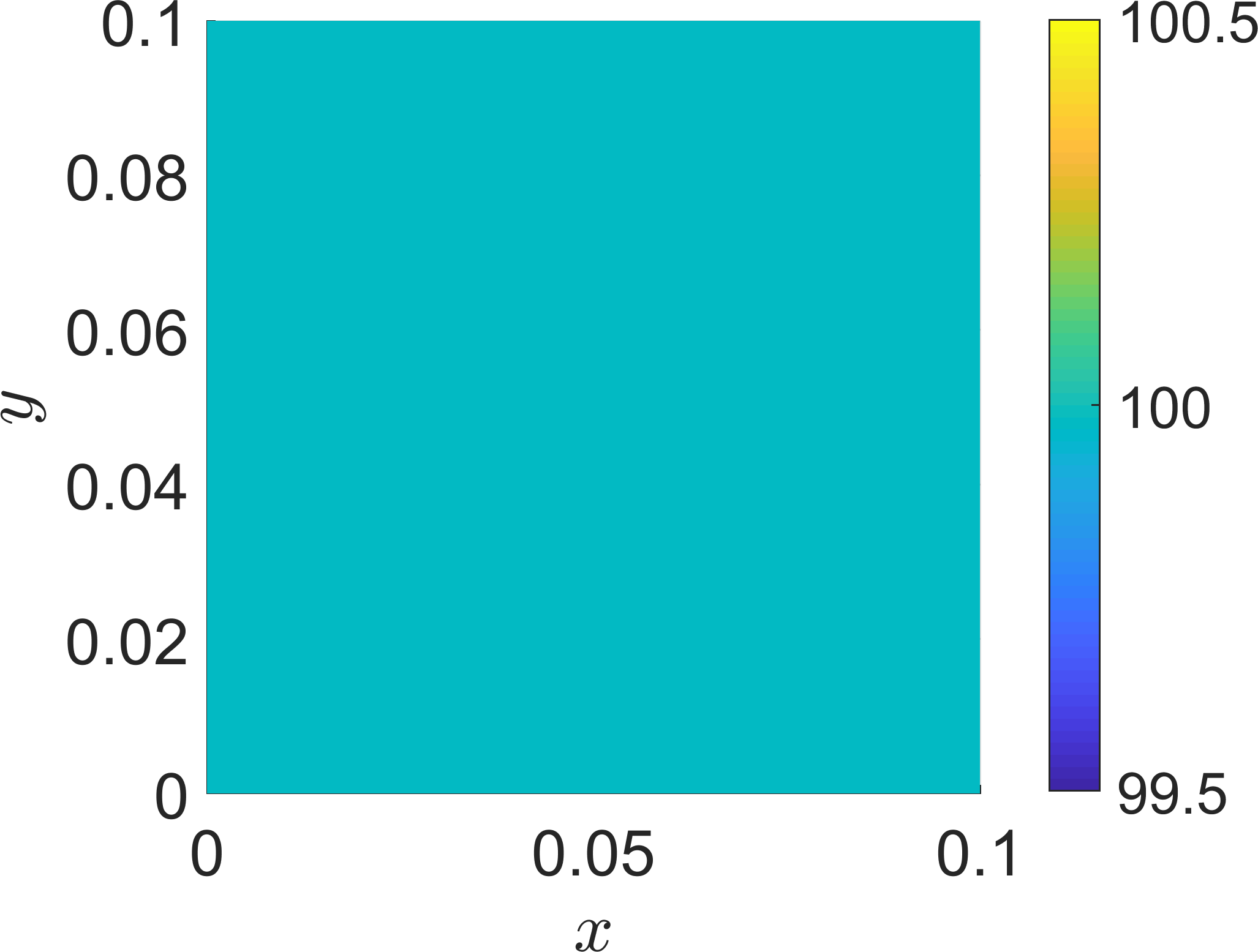}}
\subfigure[\label{hetero}]{\includegraphics[width=.4\textwidth]{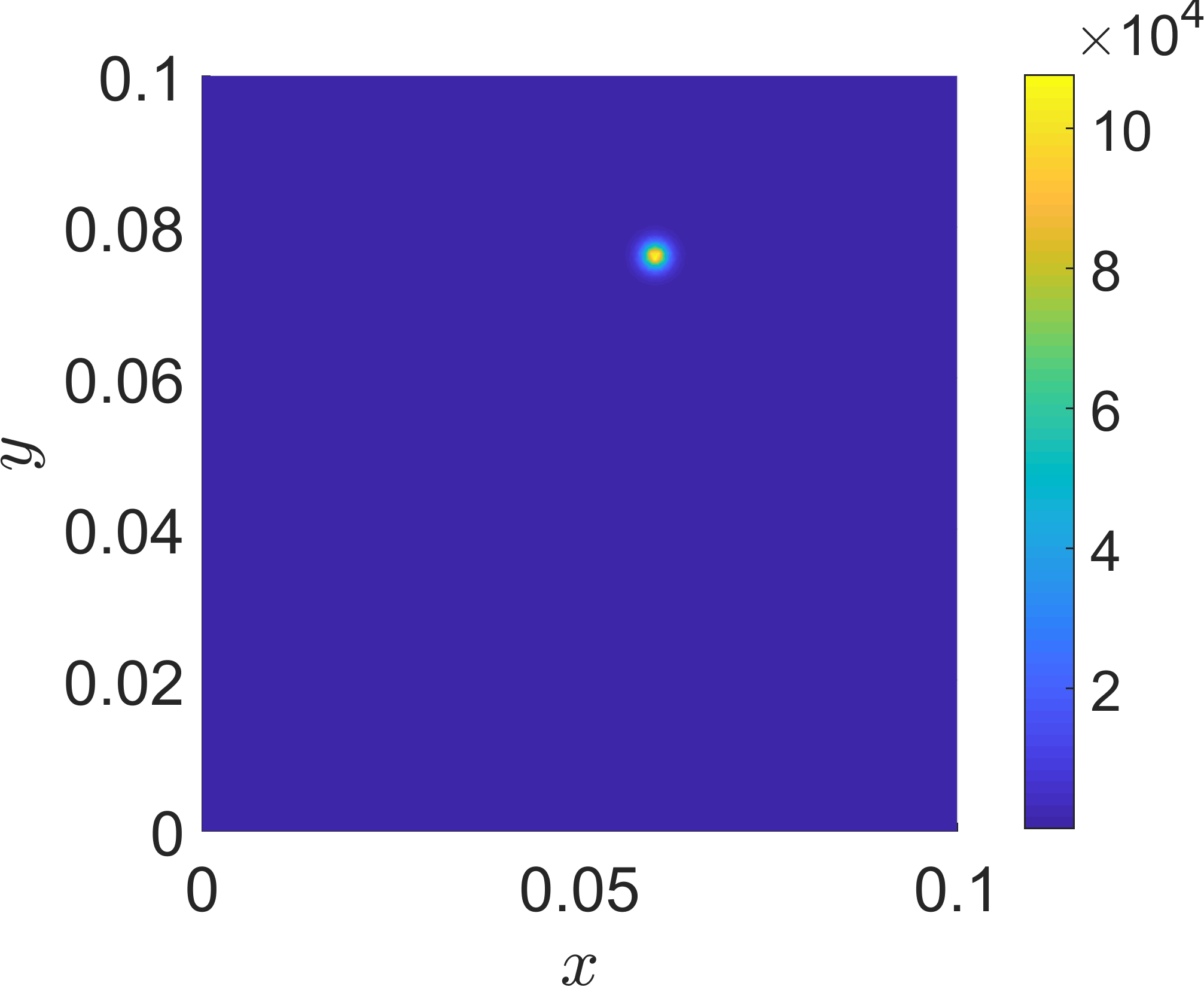}}
\caption{\label{chi_change} Simulating system \eqref{problem} with $g(u)=(1+u)^{1/2}$ and (a) $\chi(v)=\chi/v$ and (b) $\chi(v)=\chi \log(v)$, showing that blow-up does not appear to occur even when $\chi(v)$ is not constant. Parameters are $\bar{u}=100$, $\chi=10^4$ and $\sigma=10$. Due to differing scales each subfigure has its own colour bar specifying the solution density. Note that in (b), by virtue of  \eqref{MaximumPricnipleRelation} Lemma \ref{LocalExistenceLemma}, the data are taken in a such way that $v>1$, so that $\chi(v)>0$ throughout the simulation. Each image is taken at $t=100$, at which point the simulations are seen to have stopped evolving.}
\end{figure}

Next, in Figure \ref{chi_change}, we show that the solutions are bounded even when $\chi(v)$ is chosen to be non-trivial. Further, for the same parameter values, changing between $\chi(v)=\chi/v$ and $\chi(v)=\chi\log(v)$ also causes the simulations to alter between a homogeneous solution and a heterogeneous solution (compare figures \ref{homo} and \ref{hetero}).

\begin{figure}[h!!!t!!!b!!!]
\centering
\subfigure[\label{Beta_sweep_g_chi}]{\includegraphics[width=.45\textwidth]{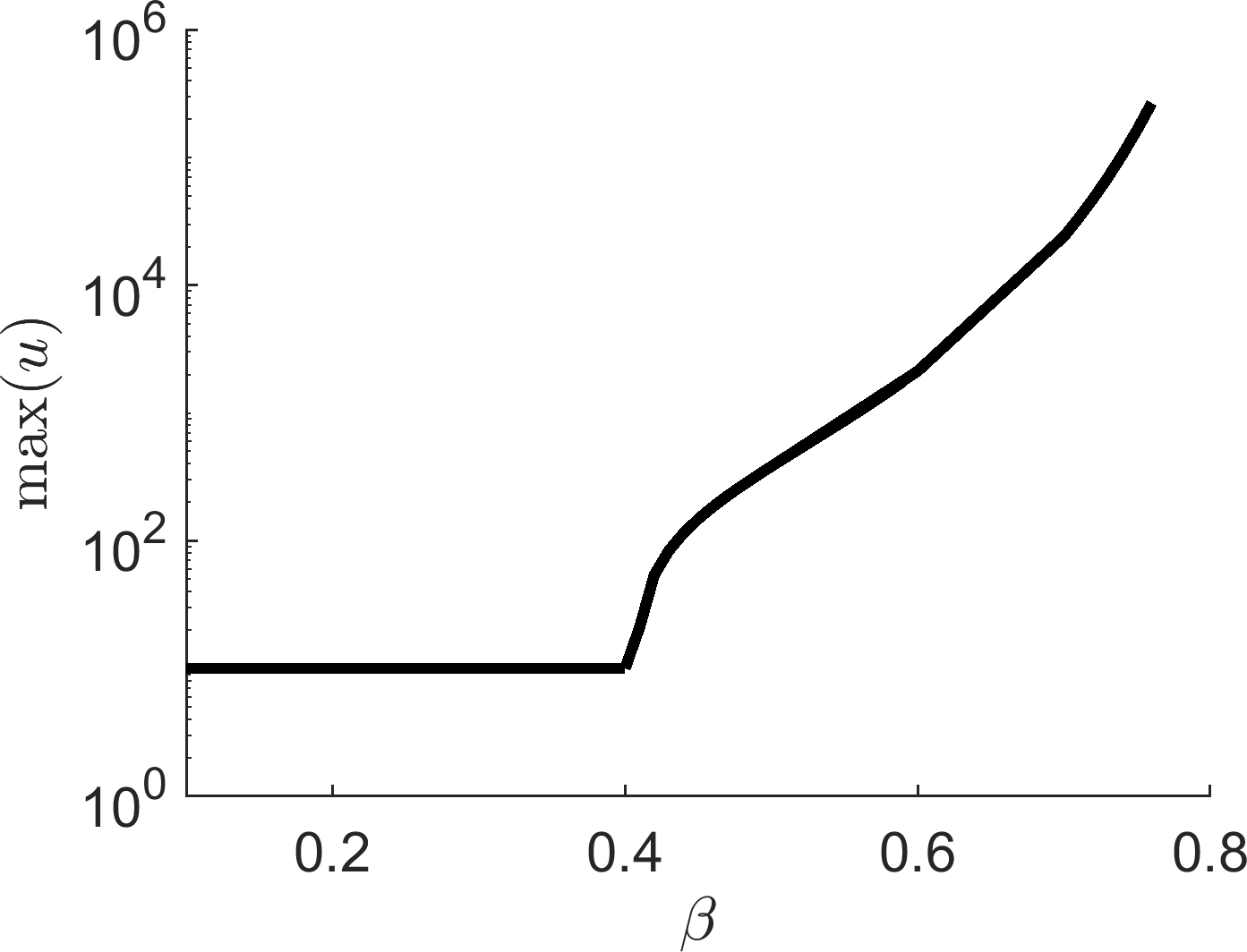}}
\subfigure[\label{Beta_sweep_g_chi_v}]{\includegraphics[width=.45\textwidth]{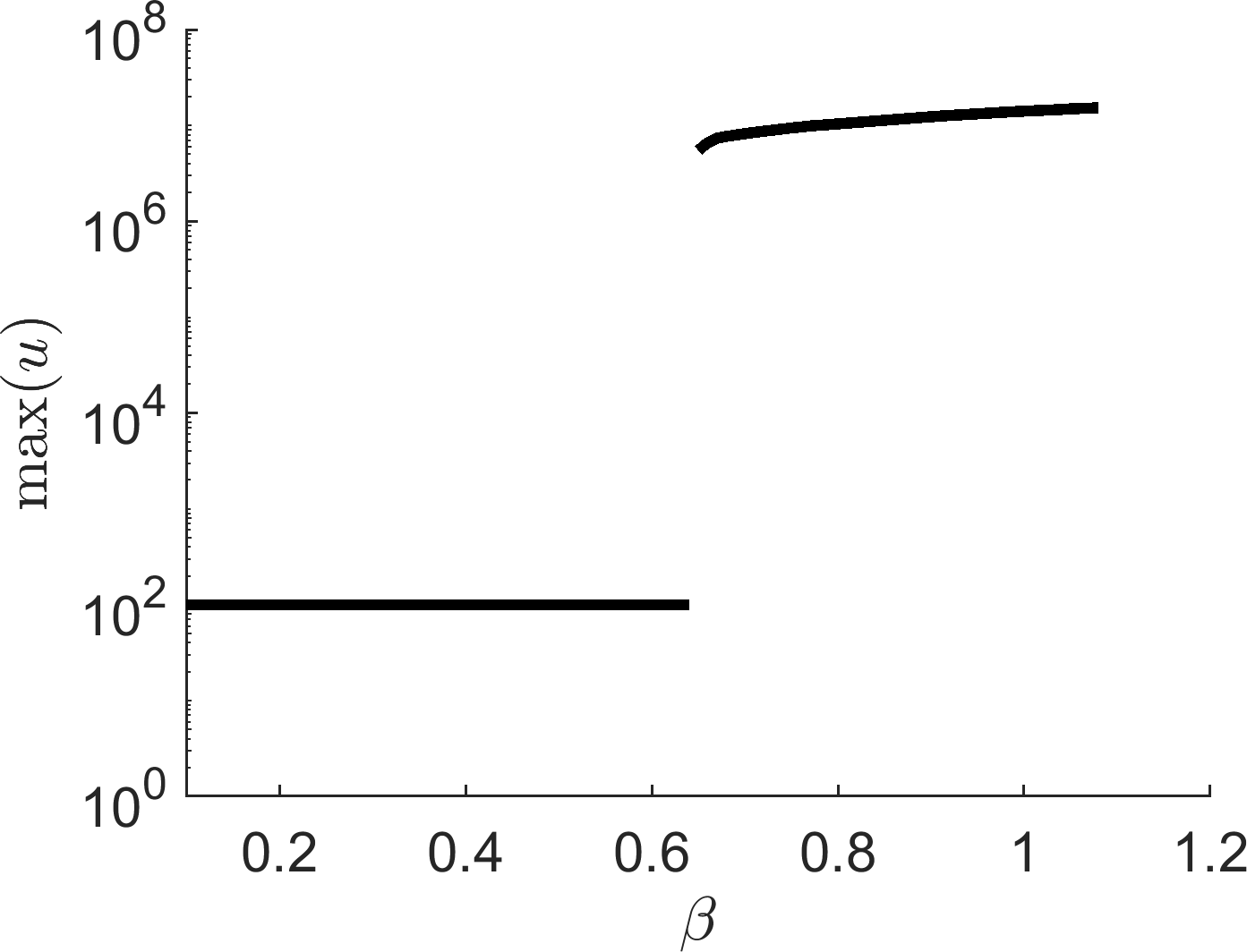}}
\caption{\label{Beta_sweep} Simulating system \eqref{problem} with $g(u)=(1+u)^{\beta}$ with increasing $\beta$ and (a) $\chi(v)=\chi=10^3$ and (b) $\chi=\chi/v=10^4/v$. Further, in (a) $\bar{u}=10$ and $\sigma=1$ and in (b) $\bar{u}=100$ and $\sigma=10$. In (a) the maximum value of $u$ appears to increase continuously with $\beta$, whilst in (b) the discontinuous jump from $\max(u)=100$ at $\beta=0.64$ to $\max(u)>10^6$ at $\beta=0.65$ illustrates that there is a bifurcation point somewhere in the interval $\beta\in[0.6,0.7]$. Each simulation was run to $t=100$, at which point the simulations had reached steady state.}
\end{figure}

Initially, and in accordance to the theoretical results, our simulations corroborated that if $g(u)=u$ the solution could become unbounded, whereas the global boundedness is guaranteed when $g(u)=(1+u)^{1/2}$. Thus, a natural question consists in analysing the possible existence of a  critical exponent, $\beta$ in $g(u)=(1+u)^{\beta}$, under which system \eqref{problem} admits stationary and homogeneous, stationary and heterogeneous, but bounded, or (even) unbounded solutions: this is addressed in the examples included in Figure \ref{Beta_sweep}. Suppose $\chi(v)=\chi>0$ and parameters for the data $\bar{u}$ and $\chi$ are chosen such that only the simulation tends to a homogeneous steady state: then a simple parameter sweep suggests that $\beta\approx 1$ is the bifurcation point (data not shown). However, if the parameters are chosen such that spike solutions exist, then the maximum values of the spikes grow as $\beta$ increases, as shown in Figure \ref{Beta_sweep_g_chi}. Here, we see that if $\beta\leq 2/5$ only the homogeneous solution is found, whereas, for $\beta\geq 2/5$ a heterogeneous solution appears. Notably, the spikes become numerically unstable for $\beta\approx 0.76$. Again, this may be due to the coarsity of the underlying mesh trying to resolve the large, but thin spikes (note that the vertical axes in Figure \ref{Beta_sweep} are logarithmic). These simulations can be compared with those seen in Figure \ref{Beta_sweep_g_chi_v}, where $\chi(v)=\chi/v$. Specifically, although we see a discontinuous jump between the homogeneous and heterogeneous solution branches at around $\beta=0.64$ (suggesting a subcritical bifurcation), we note that the maximum value of $u$ does not appear to increase much beyond $10^7$. Thus, despite the system having a spike solution, the non-trivial sensitivity produces a sort of ``controlling effect'' on the growth of the cells' distribution, $u$, even for $\beta>1$.

Over all our simulations match the analytical insights, in the ranges that are numerically feasible. Specifically, we have shown that the general properties of solutions to system \eqref{problem} heavily depend on the form of $g(u)$ and $\chi(v)$. Our numerical solutions also illustrate the need for the presented theoretical results, which confirm the existence and boundedness of solutions. Namely, since sharply peaked solutions are very difficult to numerically resolve, due to disparate scales within the system, the theoretical analysis is required to inform us as to the accuracy of our numerical schemes.

% % % % % % %
\subsubsection*{Acknowledgments}
%The author is grateful to the referees for helpful suggestions which improved this article.
GV is member of the Gruppo Nazionale per l'Analisi Matematica, la Probabilit\`a e le loro Applicazioni (GNAMPA) of the Istituto Na\-zio\-na\-le di Alta Matematica (INdAM) and is partially supported by the research project \textit{Integro-differential Equations and Non-Local Problems}, funded by Fondazione di Sardegna (2017). 

%METTERE QUELLO DEL PROGETTO YANEZ!!!!
%\bibliography{../../../Bibliography}{}
%\bibliography{Bibliography}
%\bibliography{D:/Ph.D/centralbibfile}
%\bibliography{reference}

\begin{thebibliography}{10}

\bibitem{Ascher-1997-IER}
U.~M. Ascher, S.~J. Ruuth, and R.~J. Spiteri.
\newblock {Implicit-explicit Runge-Kutta methods for time-dependent partial
  differential equations}.
\newblock {\em Appl. Numer. Math.}, 25(2-3):151--167, 1997.

\bibitem{BellomoEtAl}
N.~Bellomo, A.~Bellouquid, Y.~Tao, and M.~Winkler.
\newblock {Toward a mathematical theory of Keller--Segel models of pattern
  formation in biological tissues}.
\newblock {\em Math. Models Methods Appl. Sci.}, 25(09):1663--1763, 2015.

\bibitem{BilerLocalAndGlobal}
P.~Biler.
\newblock Local and global solvability of some parabolic systems modelling
  chemotaxis.
\newblock {\em Adv. Math. Sci. Appl.}, 8(2):715--743, 1998.

\bibitem{BrezisBook}
H.~Brezis.
\newblock {\em Functional analysis, {S}obolev spaces and partial differential
  equations}.
\newblock Universitext. Springer, New York, 2011.

\bibitem{FujieWinklerTomomiParaEllip}
K.~Fujie, M.~Winkler, and T.~Yokota.
\newblock Blow-up prevention by logistic sources in a parabolic-elliptic
  {K}eller-{S}egel system with singular sensitivity.
\newblock {\em Nonlinear Anal.}, 109:56--71, 2014.

\bibitem{FujieWinklerYokotaSignalDependent}
K.~Fujie, M.~Winkler, and T.~Yokota.
\newblock Boundedness of solutions to parabolic-elliptic {K}eller-{S}egel
  systems with signal-dependent sensitivity.
\newblock {\em Math. Methods Appl. Sci.}, 38(6):1212--1224, 2015.

\bibitem{Hillen2009UGP}
T.~Hillen and K.~J. Painter.
\newblock {A user's guide to PDE models for chemotaxis}.
\newblock {\em J. Math. Biol.}, 58(1):183--217, 2009.

\bibitem{HorstWang}
D.~Horstmann and G.~Wang.
\newblock Blow-up in a chemotaxis model without symmetry assumptions.
\newblock {\em Eur. J. Appl. Math.}, 12(2):159--177, 2001.

\bibitem{JaLu}
W.~J{\"a}ger and S.~Luckhaus.
\newblock On explosions of solutions to a system of partial differential
  equations modelling chemotaxis.
\newblock {\em Trans. Amer. Math. Soc.}, 329(2):819--824, 1992.

\bibitem{K-S-1970}
E.~F. Keller and L.~A. Segel.
\newblock Initiation of slime mold aggregation viewed as an instability.
\newblock {\em J. Theor. Biol.}, 26(3):399--415, 1970.

\bibitem{Keller-1971-TBC}
E.~F. Keller and L.~A. Segel.
\newblock {Traveling bands of chemotactic bacteria: A theoretical analysis.}
\newblock {\em J. Theor. Biol.}, 30(2):235, 1971.

\bibitem{LSUBookInequality}
O.~A. Lady\v{z}enskaja, V.~A. Solonnikov, and N.~N. Ural'ceva.
\newblock {Linear and Quasi-Linear Equations of Parabolic Type}.
\newblock In {\em {Translations of Mathematical Monographs}}, volume~23.
  American Mathematical Society, 1988.

\bibitem{Lankeit}
J.~Lankeit.
\newblock Eventual smoothness and asymptotics in a three-dimensional chemotaxis
  system with logistic source.
\newblock {\em J. Differential Equations}, 258(4):1158--1191, 2015.

\bibitem{LankeitANewApproach}
J.~Lankeit.
\newblock A new approach toward boundedness in a two-dimensional parabolic
  chemotaxis system with singular sensitivity.
\newblock {\em Math. Methods Appl. Sci.}, 39(3):394--404, 2016.

\bibitem{Nagai}
T.~Nagai.
\newblock Blowup of nonradial solutions to parabolic-elliptic systems modeling
  chemotaxis intwo-dimensional domains.
\newblock {\em J. Inequal. Appl.}, 6(1):37--55, 2001.

\bibitem{NagaiSenba}
T.~Nagai and T.~Senba.
\newblock Global existence and blow-up of radial solutions to a
  parabolic-elliptic system of chemotaxis.
\newblock {\em Adv. Math. Sci. Appl.}, 8(1):145--156, 1998.

\bibitem{NirenbergOnEllipticPDE}
L.~Nirenberg.
\newblock On elliptic partial differential equations.
\newblock {\em Ann. Scuola Norm. Sup. Pisa (3)}, 13:115--162, 1959.

\bibitem{OsYagUnidim}
K.~Osaki and A.~Yagi.
\newblock {Finite dimensional attractor for one-dimensional Keller-Segel
  equations}.
\newblock {\em Funkcial. Ekvacioj.}, 44(3):441--470, 2001.

\bibitem{SleemanLevine}
B.~D. Sleeman and H.~A. Levine.
\newblock Partial differential equations of chemotaxis and angiogenesis.
\newblock {\em Math. Methods Appl. Sci.}, 24(6):405--426, 2001.

\bibitem{OthmerStevens}
A.~Stevens and H.~G. Othmer.
\newblock Aggregation, blowup, and collapse: The abc's of taxis in reinforced
  random walks.
\newblock {\em SIAM J. Appl. Math.}, 57(4):1044--1081, 1997.

\bibitem{TelloWinkParEl}
J.~I. Tello and M.~Winkler.
\newblock A chemotaxis system with logistic source.
\newblock {\em Commun. Part. Diff. Eq.}, 32(6):849--877, 2007.

\bibitem{ViglialoroPreprintNonLinear}
G.~Viglialoro.
\newblock Global in time and bounded solutions to a parabolic-elliptic
  chemotaxis system with nonlinear diffusion and signal-dependent sensitivity.
\newblock {\em Preprint}.

\bibitem{ViglialoroVeryWeak}
G.~Viglialoro.
\newblock Very weak global solutions to a parabolic--parabolic
  chemotaxis-system with logistic source.
\newblock {\em J. Math. Anal. Appl.}, 439(1):197--212, 2016.

\bibitem{ViglialoroBoundnessVeryWeak}
G.~Viglialoro.
\newblock Boundedness properties of very weak solutions to a fully parabolic
  chemotaxis-system with logistic source.
\newblock {\em Nonlinear Anal. Real World Appl.}, 34:520--535, 2017.

\bibitem{ViglialoroWolleyDCDS}
G.~Viglialoro and T.~Woolley.
\newblock Eventual smoothness and asymptotic behaviour of solutions to a
  chemotaxis system perturbed by a logistic growth.
\newblock {\em Discrete Continuous Dyn. Syst. Ser. B.}, 22(5):1--23, 2017.

\bibitem{WinklAggre}
M.~Winkler.
\newblock Aggregation vs. global diffusive behavior in the higher-dimensional
  {K}eller-{S}egel model.
\newblock {\em J. Differential Equations}, 248(12):2889--2905, 2010.

\bibitem{WinDespiteLogistic}
M.~Winkler.
\newblock {Blow-up in a higher-dimensional chemotaxis system despite logistic
  growth restriction}.
\newblock {\em J. Math. Anal. Appl.}, 384(2):261--272, 2011.

\bibitem{WinklerFiniteTeimeBlowUpHigher}
M.~Winkler.
\newblock Finite-time blow-up in the higher-dimensional parabolic-parabolic
  {K}eller-{S}egel system.
\newblock {\em J. Math. Pures Appl.}, 100(5):748--767, 2013.

\bibitem{WinklerZAMPLogisticBlowUP}
M.~Winkler.
\newblock Finite-time blow-up in low-dimensional {K}eller--{S}egel systems with
  logistic-type superlinear degradation.
\newblock {\em Z. Angew. Math. Phys.}, 69(2):69:40, 2018.

\end{thebibliography}
%\bibliographystyle{abbrv}
%\bibliographystyle{AIMS} 
%\bibliographystyle{plainnat}

\end{document}